\documentclass[USenglish]{article}
\usepackage[utf8]{inputenc}
\usepackage[english]{babel}

\usepackage[a4paper, margin=2.5cm]{geometry}

\usepackage{microtype}

\usepackage[%
format=plain,           
labelsep=colon,               
justification=justified,      
singlelinecheck=true,         
labelfont=bf                  
]{caption}

\usepackage{subcaption}

\usepackage{amsmath}
\usepackage{amsthm}
\usepackage{amsfonts}
\usepackage{amssymb}
\usepackage{amstext}

\usepackage{mathtools}

\usepackage{physics}

\usepackage{multirow}

\renewcommand{\epsilon}{\varepsilon}

\usepackage{siunitx}
\usepackage{abstract}

\usepackage{authblk}

\title{\huge Pattern formation in clouds via Turing instabilities}
\author[1]{Juliane Rosemeier}
\author[2]{Peter Spichtinger}
\affil[1]{Institute for Atmospheric Physics, Johannes Gutenberg
    University, Mainz, Germany, rosemeie@uni-mainz.de}
\affil[2]{Institute for Atmospheric Physics, Johannes Gutenberg
    University, Mainz, Germany, spichtin@uni-mainz.de}
\date{6 October 2020}

\begin{document}

\maketitle

\renewcommand\abstractname{Abstract} 
\begin{abstract}{ Pattern formation in clouds is a well-known feature,
    which can be observed almost every day. However, the guiding
    processes for structure formation are mostly unknown, and also
    theoretical investigations of cloud patterns are quite rare. From
    many scientific disciplines the occurrence of patterns in
    non-equilibrium systems due to Turing instabilities is known,
    i.e. unstable modes grow and form spatial structures. In this
    study we investigate a generic cloud model for the possibility of
    Turing instabilities.  For this purpose, the model is extended by
    diffusion terms. We can show that for some cloud models, 
    i.e special cases of the generic model, no Turing
    instabilities are possible. However, we also present 
    a general class of cloud models, where Turing instabilities can occur.
    A key requisite is the occurrence of (weakly) nonlinear terms for accretion.
    Using
    numerical simulations for a special case of the general class of cloud models, 
    we show spatial
    patterns of clouds in one and two spatial dimensions.  
    From the numerical simulations we can see that the competition 
    between collision terms and sedimentation is an important issue 
    for the existence of pattern formation.
    }
\end{abstract}

\section{Introduction}
\label{sec:Introduction}

Pattern formation is a general feature in nature. We find patterns in
many different locations and research fields, e.g. sand ripples at
sand dunes or at the beach, stripes on zebras and fishes, convective
cells  in Rayleigh-Benard convection, spiral states in chemical
reaction systems as e.g. the famous Belousov–Zhabotinski system, 
and many other examples. The generation of structures is a common feature 
for systems out of thermodynamic equilibrium. In contrast to states at
equilibrium, which tend to be homogeneous, an external forcing driving
a system out of equilibrium has the potential to form new
structures. These structures can have different forms,
i.e. homogeneous or inhomogeneous in space and stationary or
oscillatory in time (see, e.g.,\cite{cross_hohenberg1993}). Pattern
formation is an emergent process, and is usually not  predictable a
priori from the underlying micro states of the system; the structures on 
larger scales often appear in a spontaneous way.
Research on pattern formation is an important field in many disciplines in natural 
sciences e.g. mathematical biology
\cite{murray2003}, chemistry \cite{kondepudi_prigogine2015},
fluid dynamics 
(e.g. Rayleigh-Benard convection, see \cite{bodenschatz_etal2000})
and many other fields.

There are several approaches to represent pattern formation in
models. One of the first approaches was presented by
\cite{turing1952} in his seminal article on morphogenesis. Chemical
reactions are represented by a system of ordinary differential
equations (ODEs). This set of equations is extended by diffusion
terms, i.e. a Laplacian in spatial directions is added to each
equation representing the concentration of a chemical species. It can
be shown by linear stability analysis that under certain conditions
(e.g. different diffusion coefficients) stable stationary points of the ODE
system can be destabilised, i.e. some Fourier modes become unstable
and grow, until they become saturated by nonlinear effects.  Since
only wave numbers out of a finite interval become unstable, spatial
structures become visible. This phenomenon is called Turing
instability. There are other attempts to represent structures in
models; a whole zoo of structure equations is available
 \cite{cross_hohenberg1993}. However, these approaches
are often empirical and the variables are not directly linked to
physical quantities. Sometimes, it is possible to reduce or
reformulate an underlying physical system of equations to a known
structure equation (\cite{monroy_naumis2020}). The approach of using
reaction-diffusion equations is more direct, but often ignores other
feedback due to the simplistic starting point. Nevertheless,
reaction-diffusion equations provide an important class of equations for pattern
formation, and are directly linked to the physical variables.

In atmospheric physics, a very prominent example of emerging structures 
is pattern formation
in clouds, which can be seen nicely from surface observation as well
as obtained by remote sensing techniques (e.g. from
satellites). Surprisingly, the investigation of pattern formation in
clouds is currently not a widespread topic in atmospheric physics. 
During the 1980s and 1990s several investigations and empirical studies 
on pattern formation in liquid clouds were carried out, see e.g. the review 
on cloud streets in the planetary boundary layer \cite{etling_brown1993} or 
the series on cloud clustering 
\cite{weger_etal1992, zhu_etal1992, weger_etal1993, lee_etal1994,
  nair_etal1998}.  
There are only few newer studies on pattern formation, mainly in 
connection with
investigations of open and closed cells in marine stratocumulus
(see, e.g., \cite{glassmeier_feingold2017, khouider_bihlo2019}); however, rigorous
and theoretic investigations on the formation of patterns for clouds
are lacking. 
This is surprising, since
internal structures of clouds constitute a serious uncertainty in
terms of radiative feedback. Radiative transfer in homogeneous media
is completely different than in inhomogeneous media. For the
investigations of Earth’s energy budget, clouds play a major role due
to scattering and reflection of sunlight as well as trapping infrared
radiation by absorption and re-emission. In structured clouds, many
assumptions of radiative transfer in homogeneous media do not work
anymore; for instance multiple scattering occurs frequently, and
horizontal radiative transport becomes more important. Thus, in this
respect the investigation of structured (i.e. inhomogeneous) clouds
and their origin and evolution is quite essential for meaningful
estimations of cloud radiative forcings.

There is another difficulty concerning the representation of cloud
patterns in models. Clouds constitute an ensemble of many water
particles. In cloud physics, 
one often considers processes on the scale of individual particles,
which are only partly understood until now. The description of the
statistical ensemble of cloud particles, forming the macroscopic
``object'' cloud, is not very precise, and is lacking a 
rigorous  formulation.  There are some attempts based
on Boltzmann-type evolution equations  (see,
e.g., \cite{morrison_etal2020}, \cite{beheng2010}) ,however there is no general theory of
clouds and no basic set of equations as a common ground to start is
available in cloud physics. 
In contrast when the motion of dry air shall be described, 
the Navier-Stokes equations can be used.
For the description of
clouds, often averaged variables such as the number concentration or mass
concentration of particles are used.  It is possible to relate these
quantities to general moments of the underlying size/mass distribution
of the particle ensemble  ( see, e.g. \cite{beheng2010},
\cite{khain_etal2015}).  For these averaged variables,
the process rates of the cloud processes are often formulated by
nonlinear terms.  With these parameterisation at hand the temporal
evolution of the averaged quantities can be described by  system of
ordinary differential equations, where the process rates form the
right hand side of the ODE system, also called the reaction term. 
Since a basic theory is lacking,   
the formulation of the process rates differ among the available cloud models, and
often they are not mathematically consistent.  For
instance, the uniqueness of solutions of the ODE system is not always
guaranteed, and 
often requires a more rigorous treatment (see,
e.g.,\cite{hanke_porz2020}). Nevertheless, these cloud models are often
used, and they are useful for scientific investigations as well as
for operational weather forecasts or climate predictions.

In this study, we  investigate the potential of generic cloud
models, formulated in a former study (\cite{rosemeier_etal2018}), 
to form spatial structures. We couple the model
equations with diffusion terms, i.e. Laplacians in the spatial directions;
this leads to reaction-diffusion equations for cloud physics schemes,
which will be investigated in terms of Turing instabilities.

The study is structured as follows: In the next section we will
briefly describe the generic cloud model and the represented
processes. In section~\ref{sec:LinearStabilityAnalysis} we present the
approach of linear stability theory, leading to conditions for Turing
instabilities in reaction-diffusion equations. 
The generic cloud model always allows a trivial equilibrium 
(no clouds, only rain); in section~\ref{sec:TrivialEqu} we show that
this equilibrium state cannot form Turing instabilities. 
In section~\ref{sec:noDestab} we present a special case of a cloud model,
which does not show pattern formation; this case contains standard
cloud models. In contrast, in
section~\ref{sec:CloudSchemePatternFormation} we present a general class of 
cloud models allowing Turing instabilities and thus pattern formation. 
In addition, a special case is explored for investigating several 
features of the model. 
In the  following section~\ref{sec:NumericSims} this cloud
model is numerically simulated for a 1D and a 2D scenario, and 
these results are shown. We end the study with a summary and some conclusions.

\section{Generic cloud model}
\label{sec:GenericCloudModel}

We present the generic cloud model formulated in the former study
by  \cite{rosemeier_etal2018}. The model represents clouds consisting
exclusively of liquid water droplets, so-called warm clouds.  The
droplet population is divided into two different regimes, namely cloud
droplets and rain drops, respectively. Cloud droplets are water
droplets of small sizes (radius smaller than
$\sim \SI{40}{\micro\metre}$), whereas rain drops are larger water
particles. This separation can be seen in detailed simulations
(see, e.g., \cite{beheng2010} his figure~4), and partly in
measurements.  All water particles fall in vertical direction due to
gravity. Since for small droplets the fall velocities are
very small due to friction of air, we can assume that these droplets
are stationary, in contrast to large rain drops, which fall out
faster. This separation was first proposed by  \cite{kessler1969} in
an \emph{ad hoc} manner; however, it could be justified by the
simulations mentioned above.  We consider the mass concentrations of
these two populations as variables in the model. The phase transitions
(water vapour vs. liquid water) are guided by the saturation ratio
$S_l=\frac{p_v}{p_s(T)}$, i.e. the ratio of partial pressure of water
vapour, $p_v$, and its temperature dependent saturation vapour
pressure, $p_s(T)$. Thermodynamic equilibrium, i.e. coexistence of
gaseous and liquid water, is then fulfilled at $S_l=1$; for
simplification of the notation, we also introduce the supersaturation
$S:=S_l-1$, i.e. equilibrium is reached for $S=0$.  For liquid clouds,
the following processes must be taken into account.
\begin{itemize}
\item \emph{Condensation and diffusional growth/evaporation}\\
  Cloud droplets are formed at thermodynamic conditions slightly
  beyond thermodynamic equilibrium, i.e. at supersaturation ($S>0$);
  actually, aerosol particles are activated and after passing a
  critical size, as given by Köhler theory (\cite{koehler1936}), they
  constitute cloud droplets.  In simple cloud models, this process of
  condensation is simplified and represented together with diffusional
  growth.  Cloud droplets can grow or shrink by uptake or evaporation
  of water vapour, which is provided by diffusion; this process is
  also driven by the supersaturation, which controls the thermodynamic
  equilibrium. Diffusional growth is quite inefficient for large
  droplets, thus this process is only relevant for small droplets,
  i.e. for the cloud droplet category. Both processes, condensation
  and diffusional growth (or on the contrary evaporation for $S<0$) are
  represented by a rate $C=c'Sq_c$ with a suitable constant $c'$
  depending on temperature and pressure only. For simplification, we
  assume in our investigations a permanent source of supersaturation,
  e.g. driven by a vertical upward motion. Thus, we can also neglect
  evaporation of water droplets for the system.
  \smallskip 
\item \emph{Collision processes}\\
  Since water particles fall with different velocities  depending on
  their masses, there will be collisions between neighbouring particles
  of different size, and these particles will eventually form a single 
  droplet after collision
  (so-called collision-coalescence).  Because of the artificial
  splitting of the whole droplet ensemble into two categories, we have to
  consider two (artificial) processes:\\
  (1) Two cloud droplets collide and form a large rain drop; this
  process is called autoconversion. (2) A large rain drop collects a
  small cloud droplet by collision; this process is called
  accretion. These processes are usually modelled in the spirit of
  population dynamics, using nonlinear terms; however, derivations
  from integrals over size/mass distributions lead to similar descriptions
  (see, e.g., \cite{seifert_beheng2006a}, \cite{beheng2010}. Autoconversion can be
  represented by terms $A_1=a_1q_c^\gamma$ with a suitable constant
  $a_1>0$ and an exponent $\gamma>0$.  For accretion, the terms can be
  formulated as $A_2=a_2q_c^{\beta_c}q_r^{\beta_r}$ with a suitable
  constant $a_2>0$ and exponents $\beta_c,\beta_r>0$, mimicking a
  generalised predator-prey process.
\smallskip     
\item \emph{Sedimentation of particles}\\
 For the representation of rain drops falling out of a cloud level,
  in general we would have to consider a hyperbolic term in the
  vertical direction.  For simplicity, we assume just one
  atmospheric layer with a prescribed vertical extension. Thus, we
  discretise the hyperbolic term and assume a constant flux of mass
  from above. Then, the sedimentation term can be approximated by
  $D=B-dq_r^\zeta$, with constants $B,d>0$ and an exponent $\zeta>0$.
  Note, that terminal velocities of cloud particles can be
  parameterised by power laws (see, e.g., \cite{seifert_etal2014}).
\end{itemize}
Using the representation of the processes as stated above, we obtain
the generic cloud scheme  described in  \cite{rosemeier_etal2018}:

\begin{subequations}
  \label{eq:generic_generic_model}
  \begin{align}
    \dv{q_c}{t} &= c^{\prime}Sq_c - a_1q_c^{\gamma} -
                  a_2q_c^{\beta_c}q_r^{\beta_r}, \label{eq:generic_generic_model-qc}\\   
    \dv{q_r}{t} &=  \ \quad \qquad a_1q_c^{\gamma} +
                  a_2q_c^{\beta_c}q_r^{\beta_r} +
                  B-dq_r^{\zeta}.\label{eq:generic_generic_model-qr} 
  \end{align}
\end{subequations}
To simplify the notation we will write $c$ instead of $c^{\prime} S$
in the remaining of the study.  For the analysis of the equations, we
assume constant environmental conditions, i.e. constant temperature,
pressure, and supersaturation ($S>0$), respectively.  This assumption leads to
an idealised situation, however it could be shown that similar
conditions can be encountered in the atmosphere (see,
e.g., \cite{korolev_mazin2003})for quite long times.  Assuming these
constant conditions allows us to investigate the asymptotic states of
the system.

The ODE system \eqref{eq:generic_generic_model} was discussed in
detail in \cite{rosemeier_etal2018}. In the presented work
the equations \eqref{eq:generic_generic_model}  are extended by
diffusion terms, so we obtain the following system

\begin{subequations}
  \label{eq:System}
  \begin{align}
    \dv{q_c}{t} &= c q_c -a_1 q_c^{\gamma} -a_2 q_c^{\beta_c}
               q_r^{\beta_r} \qquad \qquad \quad +  D_1 \laplacian q_c \\
    \dv{q_r}{t}  &= \ \ \quad \quad a_1 q_c^{\gamma} + a_2 q_c^{\beta_c} q_r^{\beta_r} - d
               q_r ^{\zeta} +B + D_2 \laplacian q_r .
  \end{align}
\end{subequations}
This is a reaction-diffusion system (or Turing system).  Note, that
the added diffusion terms do \emph{not} represent molecular dynamics,
as in chemical systems. Actually, these terms can be seen as a
representation of unresolved (dynamical) processes, as e.g. small
eddies or turbulence.  For the representation of turbulence in subgrid
scale schemes or entrainment due to unresolved eddies, often gradient
terms are used  (see, e.g.,\cite{deardorff1972}, \cite{stull1988}). 
This approach leads to diffusion terms in the equations for the mean variables. The
different values of the diffusion constants for the two water species
can be motivated as follows:

Small clouds droplets will mainly follow the small scale motions in the system, 
thus the diffusion coefficient $D_1$ for this species should be large. 
On the other hand, rain drops are mostly accelerated by gravity, 
thus they are less affected by small scale motions. For this species, 
the diffusion coefficient $D_2$ can be chosen different from the 
coefficient $D_1$, e.g. we would assume $D_2<D_1$. 

In the sequel the system \eqref{eq:System} is investigated with
respect to pattern formation.  The occurrence of patterns cannot be
guaranteed for the generic model, i.e. for all possible choices of
parameters, but in some cases linear stability analysis predicts
pattern formation.  These findings can be confirmed by numerical
simulations. In addition, numerical simulations 
with an extended 
parameter range might lead to further insights into potential pattern formation.

\section{Linear stability analysis}
\label{sec:LinearStabilityAnalysis}

The ideas of linear stability analysis 
(e.g. \cite{turing1952}) can be used for the determination of stable and unstable 
modes of the system of equations. A classical example for the analysis of 
reaction-diffusion equations using linear stability analysis is the investigation
of the Brusselator as a simple system describing chemical reactions
(see, e.g., \cite{cross_greenside2009} pp. 105-108).
In this section we mostly follow the exposition given in 
\cite{cross_greenside2009} for a 2D system of reaction-diffusion
equations, as, e.g., given by \eqref{eq:System}.\\
The subsequent 2D reaction-diffusion system is given by

\begin{subequations}
  \label{eq:Turing}
  \begin{align}
    \dv{u_1}{t} &= f_1 \qty(u_1,u_2) + D_1 \laplacian u_1 \\
    \dv{u_2}{t} &= f_2 \qty(u_1,u_2) + D_2 \laplacian u_2. 
  \end{align}
\end{subequations}
In a first step, we determine the stationary and homogeneous equilibrium states, 
thus we omit the diffusion terms. By neglecting the Laplacians, we
obtain a system of  
ordinary differential equations

\begin{subequations}
  \label{eq:TuringOhneDiff}
  \begin{align}
    \dv{u_1}{t} &= f_1 \qty(u_1,u_2)  \\
    \dv{u_2}{t} &= f_2 \qty(u_1,u_2) .
  \end{align}
\end{subequations}
The right hand side is called the reaction term. We want to derive
conditions for a stable equilibrium of \eqref{eq:TuringOhneDiff} which
can be destabilised by diffusion terms. First, we consider an
equilibrium solution $u_{e1},u_{e2}$ of the system
\eqref{eq:TuringOhneDiff}. By definition it satisfies the equations

\begin{subequations}
  \label{eq:TuringOhneDiff_GG}
\begin{align}
    0 &= f_1 \qty(u_{e1},u_{e2})  \\
    0 &= f_2 \qty(u_{e1},u_{e2}).
  \end{align}
\end{subequations} 
Next we compute the Jacobian of \eqref{eq:TuringOhneDiff} evaluated at
the equilibrium solution $u_{e1},u_{e2}$

\begin{align}
  \label{eq:Jacobi}
  Df |_{\qty(u_{e1}, u_{e2})} = 
  \begin{pmatrix} 
    a_{11} & a_{12} \\ a_{21} & a_{22}  \\  
  \end{pmatrix},
\end{align}
where the entries of the matrix are determined by

\begin{equation}
  \label{eq:Jacobi_partial}
  a_{ij}= \pdv{f_i}{u_j} \qty(u_{e1},u_{e2}).
\end{equation}
The (potentially complex) eigenvalues of the Jacobian at the
equilibrium state are denoted by $\sigma_1,\sigma_2 $.  The equilibrium
solution $u_{e1},u_{e2}$ is asymptotically stable if and only if
the following relations are fulfilled

\begin{subequations}
  \begin{align}
    \text{tr}\qty(Df) &\coloneqq a_{11} + a_{22} =\sigma_1+\sigma_2 < 0 \label{eq:stableEq1}\\
    \text{det}\qty(Df) &\coloneqq a_{11}  a_{22} - a_{12}  a_{21} =\sigma_1\cdot\sigma_2> 0
                         . \label{eq:stableEq2}
  \end{align}
\label{eq:stableEq}
\end{subequations}
This is equivalent to the more common condition for asymptotic
stability, i.e. $\Re(\sigma_i)<0$  for $i=1,2$.

Now we consider the system \eqref{eq:Turing} including the diffusion
terms.  For this purpose, we use spatial coordinates
$x=(x_1,\ldots,x_n)^T$ and a generalised wave number vector
$k=(k_1,\ldots,k_n)^T$. For each spatial direction, we consider linear
waves with wave lengths $\lambda_i=\frac{2\pi}{k_i}$. The Laplacian is
defined by $\displaystyle\laplacian =\sum_{i=1}^n\pdv[2]{~}{x_i}$. The
spatial dimension is given by $n\ge 1$.

\begin{figure}
    \centering
    \includegraphics[width=0.65\textwidth]{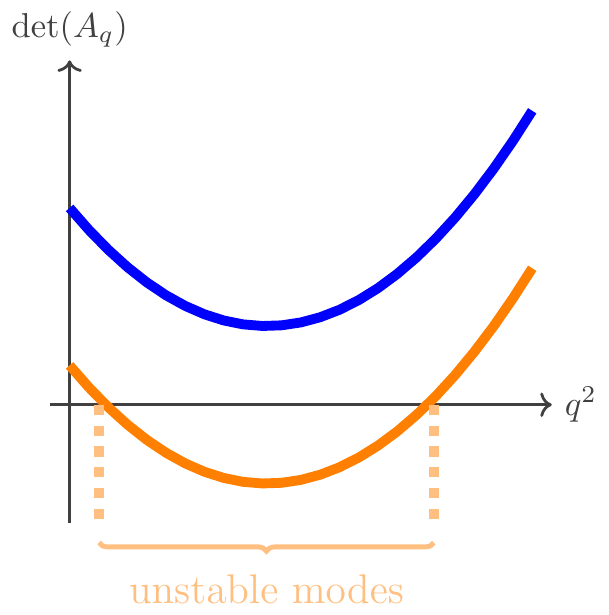}
    \caption{The quadratic polynomial \eqref{eq:quadrPloy}, which
      gives the determinant of the Jacobian
      \eqref{eq:JacobiMitDiffusion}, 
    is shown. No
  unstable modes occur if the minimum of $p_2$ is positive (blue
  line). Unstable modes are possible if the minimum of $p_2$ is
  negative (orange line).}
    \label{fig:modes}
\end{figure}

For the linear stability analysis of the reaction-diffusion system, we 
replace the reaction term by its linearisation evaluated at
$u_{e1},u_{e2}$, i.e. with the linearisation $u=u_e+u_p$, with a small 
perturbation $u_p$
around the constant equilibrium state $u_e$, we obtain

\begin{subequations}
  \label{eq:Turing_linearisiert}
  \begin{align}
    \pdv{u_{p1}}{t}&= a_{11} u_{p1} + a_{12} u_{p2} + D_1 \laplacian u_{p1}  \\
    \pdv{u_{p2}}{t}&= a_{21} u_{p1} + a_{22} u_{p2} + D_2 \laplacian u_{p2}  .
  \end{align}
\end{subequations} 
We want to derive conditions for the destabilisation of
$u_{e1},u_{e2}$ due to the diffusion terms. For simplification we
assume periodic boundary conditions; therefore a Fourier
discretisation in space with a superposition of linear wave modes
$\exp(i k x)$ can be applied.  The system
\eqref{eq:Turing_linearisiert} shall be solved by a separation
ansatz

\begin{equation}
  \label{eq:ansatz_fuer_linearisiertes_Turingmodell}
  u_p=  \begin{pmatrix} 
    u_{1q}  \\  u_{2q}  \\  
  \end{pmatrix} \exp(\sigma_q t) \exp(i k x) .
\end{equation}
using the eigenvalues $\sigma_q$ representing a single Fourier mode. 
Such a Fourier mode is an eigenfunction of the Laplacian, hence the equation

\begin{equation}
\label{eq:laplace_eigenfunction}
   \laplacian \exp(i k x) = -\left(\sum_{i=1}^nk_i^2\right) \exp(i k x) =-q^2\exp(i k x)
\end{equation}
holds, with the sum over all squared wave numbers $q^2=\sum k_i^2$, which is
used as an index.  Substituting \eqref{eq:laplace_eigenfunction} into the
linearized equation \eqref{eq:Turing_linearisiert} 
leads to the eigenvalue problem

\begin{equation}
  \label{eq:Turing_EWProblem}
  Df_q\, u_q = \sigma_q u_q ,
\end{equation}
for the coefficient $u_q$ of the Fourier mode
\eqref{eq:ansatz_fuer_linearisiertes_Turingmodell}
where the matrix $Df_q$ is given by

\begin{align}
  \label{eq:JacobiMitDiffusion}
  Df_q =
  \begin{pmatrix} 
    a_{11}-D_1 q^2 & a_{12} \\ a_{21} & a_{22} -D_2 q^2 \\  
  \end{pmatrix} .
\end{align}
For the determination of the eigenvalues $\sigma_q$ of the matrix $Df_q$
the roots of the quadratic polynomial

\begin{equation}
  \label{eq:Turing_charaktisticPolynomial}
  0 = \det(Df_q-\sigma_q I) = \sigma_q^2 - \qty(\tr(Df_q)) \sigma_q +
  \det(Df_q) .
\end{equation}
must be determined. The eigenvalues $\sigma_{qi}$ are given by

\begin{equation}
  \label{eq:Turing_eigenwerte }
  \sigma_{q1/2} = \frac{1}{2} \tr(Df_q) \pm \frac{1}{2} \sqrt{
    \qty(\tr(Df_q))^2 - 4 \det(Df_q)} .
\end{equation}
It follows that the mode which belongs to the wave number $q$ is
asymptotically stable if and only if

\begin{subequations}
\label{eq:stableEq_TuringSystem}
  \begin{align}
   \tr(Df_q) &\coloneqq  a_{11} + a_{22} - \qty(D_1 +D_2) q^2 = \sigma_{q1}+\sigma_{q2}< 0 \label{eq:stableEq_TuringSystem_tr}\\
    \det(Df_q) &\coloneqq \qty( a_{11}- D_1 q^2) \qty( a_{22}- D_2 q^2) - a_{12}  a_{21} = \sigma_{q1}\cdot\sigma_{q2}
    > 0 \label{eq:stableEq_TuringSystem_det}.
  \end{align}
\end{subequations}
Remember, that the condition can also be formulated in terms of
determinant and trace of the original ODE system, i.e.

\begin{subequations}
  \begin{align}
    \tr(Df_q) &\coloneqq
                \tr(Df)-(D_1+D_2)q^2\\
    \det(Df_q) &\coloneqq
                 \det(Df)-\qty(D_1a_{22}+D_2a_{11})q^2+D_1D_2q^4
  \end{align}
\end{subequations}
We are interested in conditions for the destabilisation of a mode $q$.
Equation \eqref{eq:stableEq_TuringSystem_tr} is satisfied because
equation \eqref{eq:stableEq1} is valid and $D_1,D_2,q^2>0$. 
The only way for destabilisation is to violate
condition \eqref{eq:stableEq_TuringSystem_det}, thus we look for a mode
which fulfills 

\begin{equation}
\label{eq:destabCond}
\qty( a_{11}- D_1 q^2) \qty( a_{22}- D_2 q^2) - a_{12}  a_{21} < 0.
\end{equation}
The left hand side of \eqref{eq:destabCond} defines a
quadratic polynomial in $q^2$,

\begin{equation}
\label{eq:quadrPloy}
p_2 (q^2) = \qty( a_{11}- D_1 q^2) \qty( a_{22}- D_2 q^2) - a_{12}a_{21} .
\end{equation}
All modes $q$ with $p_2(q^2)<0$ are unstable. The quadratic polynomial
admits a minimum at

\begin{equation}
\label{eq:Minimum_Turing}
q_m^2 = \frac{D_1 a_{22} + D_2 a_{11}}{2D_1D_2},
\end{equation} 
which constitutes the ``most unstable'' Fourier mode.
Inserting the relation \eqref{eq:Minimum_Turing} into
\eqref{eq:destabCond} yields the condition

\begin{equation}
\label{eq:unstableDiffusionSystem}
D_1 a_{22} + D_2 a_{11 } > 2 \sqrt{D_1 D_2 \qty(a_{11} a_{22} - a_{12}
  a_{21})}
\end{equation}
or in the reformulated version

\begin{equation}
  D_1 a_{22} + D_2 a_{11 } > 2\sqrt{D_1D_2\det(Df)}
\end{equation}
There is a chance to find an unstable mode $q$ if
\eqref{eq:unstableDiffusionSystem} holds, see figure
\ref{fig:modes}. The conditions \eqref{eq:stableEq1} and
\eqref{eq:unstableDiffusionSystem} can be satisfied when $a_{11}$ and 
$a_{22}$ have opposite signs.

\section{The trivial equilibrium of the generic cloud model}
\label{sec:TrivialEqu}

We now show that the trivial equilibrium of the generic cloud
model~\eqref{eq:generic_generic_model} in case of a stable stationary point
never leads to Turing instabilities. We start with the generic cloud
model~\eqref{eq:generic_generic_model}, leading to the stationary  point

\begin{equation}
    q_{ce}=0,~ q_{re}=\left( \frac{B}{d}\right)^{\frac{1}{\zeta}}.
\end{equation}
Since this implies no cloud, just rain, in the atmospheric layer, this
state is called trivial equilibrium.  Actually, this stationary point is only
valid for linear stability analysis for values of the exponents
$\gamma\ge 1$ and $\beta_c\ge 1$. 
Otherwise, the partial derivatives with respect to $q_c$ do not exist at $q_c=0$. 
Cloud models with $\gamma < 1$ or $\beta_c < 1$ also lack Lipschitz
continuity. Therefore, the Picard-Lindelöff theorem does not guarantee
the unique solvability of the ODE system when the initial value for
$q_c$ is given by  $q_c=0$.  
For discussions of such cloud
models and possible extensions to uniqueness, see the recent study by
\cite{hanke_porz2020}. In the remainder of the study developed here we
will always  
assume $\gamma\ge 1,\beta_c\ge 1$. For linear stability analysis we have to
consider the Jacobian $Df|_{(q_{ce},q_{re})}$. As discussed by
\cite{rosemeier_etal2018} the Jacobian always has the form

\begin{equation}
   Df|_{(q_{ce},q_{re})}= \begin{pmatrix} 
    a_{11} & 0 \\ a_{21} & a_{22} \\  
  \end{pmatrix} \coloneqq A, 
\end{equation}
with
$a_{22}=-d\zeta \left(\frac{B}{d}\right)^{\frac{\zeta-1}{\zeta}}<0$.
As $\beta_c \ge 1 $ is assumed, we obtain $a_{12}=0$.
If even the conditions $\gamma>1$ and $\beta_c>1$ hold, we obtain $a_{21}=0$,
otherwise $a_{12}>0$; for details, see calculations in 
appendix A. 
Nevertheless, it is clear that the eigenvalues $\sigma_i$ are given by
$\sigma_1=a_{11}, \sigma_2=a_{22}<0$ and thus
$\det(A)=a_{11}a_{22}=\sigma_1\sigma_2$, and
$\Tr(A)=a_{11}+a_{22}=\sigma_1+\sigma_2$, respectively. For a stable
stationary point, both (real) eigenvalues must be negative, leading to the
criteria \eqref{eq:stableEq}; this might be fulfilled for the choice
of parameters $c<a_1$ in case of $\gamma=1$; otherwise the stationary point 
cannot be stable (see also appendix~A).  
However, the stable stationary point can not lead
to Turing instabilities via destabilisation.  The criterion for the
existence of destabilisation \eqref{eq:unstableDiffusionSystem} can be
reduced to the following form:

\begin{equation}
   D_1 a_{22} + D_2 a_{11 } > 2 \sqrt{D_1 D_2 \qty(a_{11} a_{22})}. 
\end{equation}
Since $a_{11}=\sigma_1<0$ and $a_{22}=\sigma_2<0$, this leads to a
contradiction. This proves that the trivial  stationary point (if it exists) 
cannot be destabilised by diffusion, and thus it cannot serve for Turing
instabilities.

From a physical point of view, in this situation the source for cloud
droplets represented by the term $cq_c$ is too weak and 
collision processes (terms $A_1$ and $A_2$) reduce the cloud water such
efficiently that diffusion cannot change the quality of the stable
stationary point (i.e. no cloud with rain).

If the parameters $c,a_1$ are chosen such that
$a_{11}=\sigma_1=c-a_1>0$, an unstable equilibrium state can be
obtained. Physically, this means that condensation and diffusional
growth is much stronger than autoconversion, i.e. more cloud mass is
generated by condensation than is lost by collision processes.
This situation usually occurs in scenarios with a persistent
updraught, leading to a steady source of supersaturation and thus
permanent
cloud droplet formation and growth. 
In case of an unstable trivial equilibrium, 
the modified Jacobian at the equilibrium
admits the following form:

\begin{equation}
   Df_q|_{(q_{ce},q_{re})}=  \begin{pmatrix} 
    a_{11}-D_1q^2 & 0 \\ a_{21} & a_{22}-D_2q^2 \\  
  \end{pmatrix}.
\end{equation}
The second eigenvalue $\sigma_{q2}=a_{22}-D_1q^2$ is still negative, since $D_1,q^2 \ge 0$.
The first (real) eigenvalue can be negative for Fourier modes $q$ fulfilling the following condition

\begin{equation}
    \sigma_{q2}=a_{11}-D_1q^2<0\Leftrightarrow a_{11}<D_1q^2\Leftrightarrow 
    \frac{a_{11}}{D_1}<q^2.
\end{equation}
Thus, the absolute value of the diffusion constant $D_1>0$ decides
about the stability of the modes.

\section{A case without destabilisation}
\label{sec:noDestab}

We set $\gamma = 1$ and $\beta_c = 1$ in the system
 \eqref{eq:generic_generic_model} and show that it is not 
 possible to destabilise an asymptotically stable equilibrium 
 of this model by diffusion terms with arbitrary coefficients $D_1,D_2>0$.
The cloud scheme of the operational numerical weather prediction model
COSMO (\cite{cosmo_doc_physical_parameterization}) of the 
German weather service (DWD)
and the research model by \cite{wacker1992} admit this special form of 
the cloud scheme examined in the sequel. Particularly we consider

\begin{subequations}
  \label{eq:general_case}
  \begin{align}
    \dv{q_c}{t} &= cq_c - a_1q_c - a_2q_c q_r^{\beta_r},
    \label{eq:general_case-qc}\\ 
    \dv{q_r}{t} &=  
    \qquad \ \ a_1q_c + a_2q_c q_r^{\beta_r} + B-dq_r^{\zeta}.
    \label{eq:general-case-qr}
  \end{align}
\end{subequations}
Besides the trivial equilibrium state (see discussion in
section~\ref{sec:TrivialEqu}) 
the only non-trivial equilibrium of the system \eqref{eq:general_case}
is given by (see \cite{rosemeier_etal2018})

\begin{equation}
q_{ce} = \frac{d}{c}
\qty(\frac{c-a_1}{a_2})^{\frac{\zeta}{\beta_r}} - \frac{B}{c}
\qquad q_{re} = \qty(\frac{c-a_1}{a_2})^{\frac{1}{\beta_r}}.
\label{Struk_bildung_GG}
\end{equation}
Note, for the existence of this (non-negative) 
equilibrium state two conditions must be fulfilled, i.e.
\begin{equation}
    c>a_1,\text{~and~} d\qty(\frac{c-a_1}{a_2})^{\frac{\zeta}{\beta_r}}>B.
\end{equation}
Physically, this means that, as before, the cloud droplet source $cq_c$ is stronger 
than the sink of autoconversion. Additionally, the rain flux from above $B$ 
must not be too strong, otherwise no equilibrium state is reached, 
i.e. the rain will collect almost all cloud droplets.

Again we compute the Jacobian at the equilibrium state

\begin{align}
  \label{eq:Jacobi_}
  Df |_{\qty(q_{ce}, q_{re})} =
  \begin{pmatrix} 
    a_{11} & a_{12} \\ a_{21} & a_{22}  \\  
  \end{pmatrix} ,
\end{align}
where

\begin{subequations}
\begin{align}
a_{11} &= c-a_1-a_2 q_{re}^{\beta_r} = 0,\\
a_{12} &= - a_2 \beta_r q_{ce} q_{re}^{\beta_r-1} <0, \\
a_{21} &= a_1+a_2 q_{r,e}^{\beta_r} = c >0,\\   
\begin{split}
a_{22} &= a_2 \beta_r q_{ce} q_{re}^{\beta_r -1} -d \zeta
q_{re}^{\zeta -1} \\
&= a_2 \beta_r \qty( \frac{d}{c}
\qty(\frac{c-a_1}{a_2})^{\frac{\zeta}{\beta_r}} - \frac{B}{c})
\qty(\frac{c-a_1}{a_2})^{\frac{\beta_r - 1}{\beta_r}} - d \zeta
\qty(\frac{c-a_1}{a_2})^{\frac{\zeta -1}{\beta_r}}. 
\end{split}
\end{align}
\end{subequations}
As $a_{11} =0$, condition \eqref{eq:stableEq1} gives $a_{22} < 0$; 
note, that condition \eqref{eq:stableEq2} is fulfilled. On
the other hand in this case, condition
\eqref{eq:unstableDiffusionSystem} reduces to

\begin{equation}
D_1 a_{22} > 2 \sqrt{D_1 D_2 \qty(-1) a_{12} a_{21}}.
\label{eq:Bed_destab_wacker-COSMO}
\end{equation}
This yields a contradiction. Thus, for schemes with $\gamma=\beta_c=1$
pattern formation via Turing instabilities is impossible.

Generally, it is of interest if and when the matrix entry $a_{22}$
changes its sign, since this entry determines the quality of the stationary  point.
For this purpose we further investigate $a_{22}$ in detail

\begin{subequations}
\begin{align}
\begin{split}
a_{22} &= \frac{a_2 \beta_rd}{c}
\qty(\frac{c-a_1}{a_2})^{\frac{\beta_r +\zeta-1}{\beta_r}} - d \zeta
\qty(\frac{c-a_1}{a_2})^{\frac{\zeta -1}{\beta_r}} - \frac{a_2
  \beta_r}{c}  B \qty(\frac{c-a_1}{a_2})^{\frac{\beta_r-1}{\beta_r}}
 \\
&= \qty(\frac{ \beta_rd}{c} \qty(c-a_1)) -d \zeta 
\qty(\frac{c-a_1}{a_2})^{\frac{\zeta -1}{\beta_r}}
\underbrace{- \frac{a_2\beta_r}{c}  B
  \qty(\frac{c-a_1}{a_2})^{\frac{\beta_r-1}{\beta_r}}}_{< 0} .
\end{split}
\end{align}
\end{subequations}
For the determination of the sign of $a_{22}$ the following term must
be examined

\begin{equation}
\label{eq:Bedingung_an_einen_Term}
\frac{ \beta_rd}{c} \qty(c-a_1)) -d \zeta = d \qty(\beta_r -
\frac{a_1 \beta_r}{c} - \zeta).
\end{equation}
We can now conclude that $\beta_r \le \zeta$ is sufficient for $a_{22}$ to be negative. This condition holds for the Wacker and COSMO schemes.  

The trivial equilibrium state is given by 

\begin{equation}
q_c = 0, \qquad q_r = \qty(\frac{B}{d})^{\frac{1}{\zeta}}.
\end{equation}
It is unstable for the COSMO and Wacker schemes, but several modes $q$ are
stabilised by the diffusion terms, depending on the diffusion constant $D_1$.
However it is not generally unstable for an arbitrary choice for the
prefactors. The stable case was already discussed in
section~\ref{sec:TrivialEqu}, it cannot trigger Turing instabilities.

\section{A cloud scheme with pattern formation}
\label{sec:CloudSchemePatternFormation}

In this section we present a general class of cloud schemes of the 
form \eqref{eq:generic_generic_model} which allow pattern formation 
via Turing instabilities. Again the method
described in section \eqref{sec:LinearStabilityAnalysis} is applied. 
Cloud schemes of the form

\begin{subequations}
  \begin{align}
    \dv{q_c}{t} &= c q_c - a_1 q_c - a_2 q_c^{\beta} q_r^{\beta}\\
    \dv{q_r}{t} &=\ \quad \quad \ a_1 q_c + a_2 q_c^{\beta} q_r^{\beta} -d q_r 
  \end{align}
  \label{eq:CloudScheme_WithPatterns_general}
\end{subequations}
are considered, where $\beta > 1$. 
Thus, accretion is parameterised by the term $a_2 q_c^{\beta} q_r^{\beta}$, 
which can also be found in some standard cloud models 
 (see, e.g., \cite{khairoutdinov_kogan2000} with $\beta=1.15$ as used in the
IFS).  
For simplification, we additionally assume 
a linear autoconversion ($\gamma=1$) according to \cite{kessler1969}.
However, we will discuss later that this restriction is not crucial.
Note, that we also assume $\beta=\beta_r>\zeta=1$ as indicated in 
section~\ref{sec:noDestab}. Finally, we first omit the constant 
rain flux from above (i.e. the term $B$) for simplification; we will discuss the
inclusion of this term at the end of this section and also in
section~\ref{sec:NumericSims} for a special set of parameters.
The ODE system \eqref{eq:CloudScheme_WithPatterns_general} is extended by diffusion terms and the resulting  reaction-diffusion system is given by

\begin{subequations}
  \begin{align}
 \label{eq:CloudScheme_WithPatterns__generalTuring1}
    \pdv{q_c}{t} &= c q_c - a_1 q_c - a_2 q_c^{\beta} q_r^{\beta} \quad\qquad  
    + D_1 \laplacian q_c\\
    \pdv{q_r}{t} &=\ \ \quad \quad a_1 q_c + a_2 q_c^{\beta} q_r^{\beta} -d q_r \ 
    + D_2 \laplacian q_r.
 \label{eq:CloudScheme_WithPatterns__generalTuring2}
  \end{align}
 \label{eq:CloudScheme_WithPatterns__generalTuring}
\end{subequations}
The system \eqref{eq:CloudScheme_WithPatterns_general} has the
following nontrivial  equilibrium

\begin{equation}
q_{re} = \frac{c}{d} q_{ce} \ , \qquad \qquad 
q_{ce} = \qty(\frac{d^{\beta}}{c^{\beta}} \frac{c-a_1}{a_2})^{\frac{1}{2 \beta -1}} .
\label{sec6:equi_general}
\end{equation}
To guarantee the existence of a positive equilibrium we assume 

\begin{equation}
c>a_1 . 
\label{sec6:cond1}
\end{equation}
This is the first constraint on the admissible set of parameters. 
As before, condensation is dominant over autoconversion of cloud droplets. 
The Jacobian evaluated at the above mentioned equilibrium has the  entries

\begin{subequations}
\begin{align}
a_{11} &= \qty(1- \beta) (c-a_1) < 0 \\
a_{12} &= - \beta \frac{d}{c} \qty(c-a_1) < 0 \\
a_{21} &= a_1 + \beta (c-a_1) > 0\\
a_{22} &= d \qty(\beta \frac{c-a_1}{c} -1) .
\end{align}
\end{subequations}
 Turing instabilities can arise if  $a_{22} >0 $. This is equivalent to the condition
 \begin{equation}
     1< \beta \frac{c-a_1}{c}.
     \label{sec6:cond2}
 \end{equation}
 Note that condition~\eqref{sec6:cond2} is equivalent to the formulation
 $1-\frac{1}{\beta}>\frac{a_1}{c}$, which implies $\beta>1$ as already assumed;
 thus, the prefactors $a_1, c > 0$ can be chosen accordingly.

 The trace of the Jacobian is given by
 \begin{equation}
     \tr(Df_q)=a_{11} +a_{22} = \qty(1- \beta) (c-a_1) + 
     d \qty(\beta \frac{c-a_1}{c} -1) .
 \end{equation}
 The condition \eqref{eq:stableEq1} for negative trace holds when $d$ is chosen small enough.
 It can be shown that for the system \eqref{eq:CloudScheme_WithPatterns_general} 
 the condition on the positive determinant \eqref{eq:stableEq2} is equivalent to 
 $0 < 2 \beta -1$ which is always satisfied for $\beta > 1$.
 Consequently it is possible to chose $a_1, c$ and $d$ as well as $D_1$ and $D_2$ 
 in \eqref{eq:CloudScheme_WithPatterns__generalTuring} such that
 Turing instabilities can arise. 

As a concrete example for investigating the details we consider the case $\beta=2$, 
i.e. the cloud scheme 

\begin{subequations}
  \begin{align}
    \dv{q_c}{t} &= c q_c - a_1 q_c - a_2 q_c^2 q_r^2 \label{eq:CloudScheme_WithPatterns2_qc}\\
    \dv{q_r}{t} &=\ \quad \quad a_1 q_c + a_2 q_c^2 q_r^2 -d q_r . \label{eq:CloudScheme_WithPatterns2}
  \end{align}
  \label{eq:CloudScheme_WithPatterns}
\end{subequations}

The corresponding reaction-diffusion system has the form

\begin{subequations}
  \begin{align}
    \pdv{q_c}{t} &= c q_c - a_1 q_c - a_2 q_c^2 q_r^2 \quad\qquad  + D_1 \laplacian q_c\\
    \pdv{q_r}{t} &=\ \quad \quad a_1 q_c + a_2 q_c^2 q_r^2 -d q_r \ + D_2 \laplacian q_r.
  \end{align}
 \label{eq:CloudScheme_WithPatterns_Turing}
\end{subequations}

Applying relation \eqref{sec6:equi_general} gives the non-trivial equilibrium

\begin{subequations}
  \label{eq:CloudScheme_WithPatterns_Equilibrium}
  \begin{align}
    q_{ce} &= \qty(\frac{c-a_1}{a_2})^{\frac{1}{3}}
              \qty(\frac{d}{c})^{\frac{2}{3}}\\ 
    q_{re} &= \qty(\frac{c}{a_2d})^{\frac{1}{3}} \qty(c-a_1)^{\frac{1}{3}} 
  \end{align}
\end{subequations}
for $c>a_1$ as indicated in equation~\eqref{sec6:cond1}.

In the next step the Jacobian of \eqref{eq:CloudScheme_WithPatterns} evaluated 
at the equilibrium \eqref{eq:CloudScheme_WithPatterns_Equilibrium} is
considered.  
The condition $a_{22} > 0$ is equivalent to 
\begin{equation}
  \label{eq:CloudScheme_WithPatterns_cond2}
  c> 2 a_1.
\end{equation}
This amounts to a further constraint on the prefactors.
Additionally,  the trace of the Jacobian must be negative when
\eqref{eq:stableEq1} is supposed to hold, i.e. 
\begin{equation}
\label{eq:CloudScheme_WithPatterns_TestTraceCond}
a_{11}  + a_{22}  =
\underbrace{a_1-c}_{<0} + \underbrace{ d \qty( 1- \frac{2
    a_1}{c})}_{>0} .
\end{equation}
When $d$ is chosen sufficiently small, i.e.
\begin{equation}
  \label{eq:CloudScheme_WithPatterns_cond3}
   d < c \frac{c-a_1}{c- 2a_1} ,
\end{equation}
the relation \eqref{eq:stableEq1} holds.

\begin{figure}
    \centering
    \includegraphics[width=0.75\textwidth]{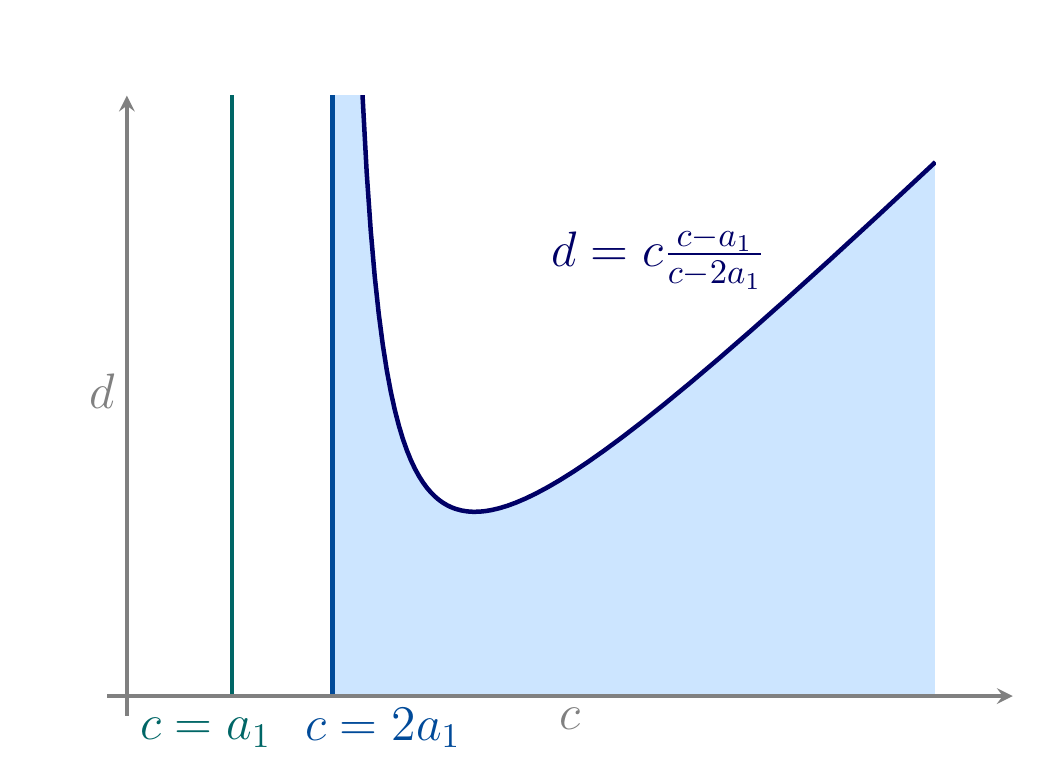}
    \caption{Constraints on the prefactors. We assume that $a_1$ is given.}
    \label{fig:constrains}
\end{figure}

Thus with an appropriate choice of $d$ according to 
\eqref{eq:CloudScheme_WithPatterns_cond3} we can satisfy the conditions
\eqref{eq:stableEq}, and  \eqref{eq:unstableDiffusionSystem}
can be fulfilled for a proper choice of the diffusion constants $D_1,D_2$. 
In summary, we have derived three limiting conditions 
\eqref{sec6:cond1}, \eqref{eq:CloudScheme_WithPatterns_cond2}, and
\eqref{eq:CloudScheme_WithPatterns_cond3},
which are
illustrated in figure \ref{fig:constrains}. For values of $c$ and $d$ in the blueish 
domain of the parameter space, the equilibrium state is stable 
and in general allows Turing instabilities.

From eq.~\eqref{eq:CloudScheme_WithPatterns_cond3} as well as from the 
phase diagram in
figure~\ref{fig:constrains} we see that the strength of the sedimentation 
(parameter $d$) 
plays a major role for the existence of Turing instabilities. 
If sedimentation is too strong
 compared to condensation, diffusion is not effective enough to distribute
the cloud spatially for generating instabilities. 

\noindent\textbf{Remarks:}
\begin{enumerate}
\item If we investigate the first equation of the generic ODE
  System~\eqref{eq:generic_generic_model-qc} 
  we can identify a relation satisfied by the non-trivial equilibrium 
  $q_{ce}, q_{re}$ which holds for any admissible set of parameters.
  If the condensation term in the original
  formulation is slightly extended to obtain the final equation
  
  \begin{equation}
    \dv{q_c}{t} = cq_c^{\gamma} - a_1q_c^{\gamma} -
                  a_2q_c^{\beta_c}q_r^{\beta_r} ,
  \end{equation}
  i.e. the condensation and the autoconversion term have the same
  exponential behaviour, 
  then we detect the following identity for the nontrivial steady states 
  satisfying $\dv{q_c}{t}=0$
  \begin{equation}
    \label{eq:conservation}
    \frac{c-a_1}{a_2}=q_{ce}^{\beta_c-\gamma}q_{re}^{\beta_r}.
  \end{equation}
  This means that a certain combination of $q_{ce} $ and $q_{re} $ is constant 
  in the  set  of the admissible exponents $\beta_c, \beta_r$ and $\gamma$, 
  provided the prefactors $a_1,a_2$ and $c$ are held constant. Therefore, 
  the value $ \frac{c-a_1}{a_2} $ will be denoted as a conserved quantity 
  in the sequel. Note that the validity of the identity \eqref{eq:conservation} 
  is unaffected by the inclusion of a term $B$ in 
  \eqref{eq:generic_generic_model-qr}.
\item In case of $\gamma=1$, we obtain for the conserved quantity in
  equation~\eqref{eq:conservation} the term
  
  \begin{equation}
    \frac{c-a_1}{a_2}=q_c^{\beta_c-1}q_r^{\beta_r} .
    \label{sec6:conserved}
  \end{equation}
  This leads to a strong simplification of the first column in the
  Jacobi matrix $Df$, since with relation \eqref{sec6:conserved} 
  the dependence on  $q_{ce}, q_{re}$ can be eliminated easily 
  ($a_{11}=\qty(1-\beta_c) \qty(c-a_1)$, $a_{21}=a_1+\beta_c(c-a_1)$). 
  This property holds for all
  values of $\beta_r>0$.
\item Including the rain flux from above $B$ into equation
 \eqref{eq:CloudScheme_WithPatterns__generalTuring2}
 leads to the modification
  
  \begin{equation}
    \dv{q_r}{t} = a_1 q_c + a_2 q_c^\beta q_r^\beta +B-d q_r .
    \label{sec6:rmk3_qr}
  \end{equation}
  Thus, the determination of the non-trivial stationary point becomes more
  complicated. For the special case $\beta=2$ we investigate \eqref{sec6:rmk3_qr} explicitely. The stationary points shall be computed.
  Following the strategy at the beginning of the section, we
  add the two ODE equations \eqref{eq:CloudScheme_WithPatterns2_qc} and \eqref{sec6:rmk3_qr}, which were set to zero before. This 
  leads to the problem of finding real roots of the
  general cubic polynomial
  
  \begin{equation}
    p_3(q_{re})=d q_{re}^3-Bq_{re}^2-c\frac{c-a_1}{a_2}.
  \end{equation}
  Using Cardano's formulas (see, e.g., \cite{bosch2018}, chapter 6),
  we can determine the roots of the polynomial directly using the
  following terms:
  
\begin{eqnarray}
  \label{eq:cardano-p}
  p & = & -\frac{1}{3}\qty(\frac{B}{d})^2\\
  \label{eq:cardano-q}
  q & = & -\frac{2B^3+27d^3c\frac{c-a_1}{a_2}}{27d^3}\\
  \label{eq:cardano-delta}
  \Delta & = & \qty(\frac{q}{2})^2+\qty(\frac{p}{3})^3
      =  \frac{27d^2c^2\frac{\qty(c-a_1)^2}{a_2^2}-4B^3c\frac{c-a_1}{a_2}}{108d^4}
\end{eqnarray}
whereas the parameter $\Delta$ decides about the quality of the roots
(e.g. one real root and two complex conjugates for $\Delta >0$). One
real root is given by

\begin{eqnarray}
  u =  \sqrt[3]{-\frac{q}{2}+\sqrt{\Delta}}, &
  v  =  \sqrt[3]{-\frac{q}{2}-\sqrt{\Delta}}, & 
  q_{re}=u+v+\frac{B}{3d}.
  \label{sec6:real_root}
\end{eqnarray}
Since the sign of parameter $\Delta$ decides about the number of real
roots, there is in general a bifurcation at $B_1$, which can be
calculated using equation~\eqref{eq:cardano-delta}:

\begin{equation}
  27d^2c^2\frac{\qty(c-a_1)^2}{a_2^2}=4B_1^3c\frac{c-a_1}{a_2}
  \Leftrightarrow
  B_1=\sqrt[3]{\frac{27}{4}d^2\frac{c(c-a_1)}{a_2}}.
\end{equation}
However, the condition~\eqref{eq:unstableDiffusionSystem} for the
existence of a Turing instability might be violated at different
values of $B$ 
as can be seen in the numerical simulations in the next section. The
equilibrium states can be inserted into the Jacobi matrix for
determining the eigenvalues. 
Using the relation \eqref{sec6:conserved}, we see that
the entries $a_{11}$ and $a_{21}$ do not depend on
$B$. Therefore when the impact of $B$ on the existence
of Turing instabilities shall be investigated, we only
have to
consider the entries $a_{12}(B),a_{22}(B)$. The sign
of $a_{22}(B)$ 
is again the key
parameter, and depends on the rain flux $B$.  Actually, for a certain
setting of parameters $c,a_1,a_2,d$ we will determine the qualitative
behaviour and the possibility of Turing instabilities numerically (see
next section).
\item The choice of exponents $\beta_c>1,\beta_r>1$  in the accretion term is motivated 
by  existing models as e.g. the operational weather forecasting model
IFS (\cite{ifs_doc_physical_parameterization}), which contains such exponents.
However, since the representation of collision processes in bulk models
is not well-defined from a basic theory, there is actually no restriction
of the choice of parameters. It is conceivable that the parameters also vary for different 
environmental regimes.
\item A nonlinear autoconversion generally affects the nontrivial stationary point as well
as the entries $a_{12},a_{21}$ of the Jacobian. However, an analytical derivation of 
general conditions for the occurrence of Turing
instabilities 
appears at least cumbersome through the nonlinear equation for the equilibrium, i.e. 
\begin{equation}
    d-a_1\qty(\frac{d}{c})^{\gamma}q_r^{\gamma-1}
    -a_2\qty(\frac{d}{c})^{\beta}q_r^{2\beta-1}=0,~~\beta>1,\gamma>1.
\end{equation}
Numerical studies for different values of parameters $\beta,\gamma$ indicate that pattern
formation is not restricted to the case $\gamma =1$. 
We could detect patterns in two simulations where $\beta_c = \beta_r = 1.25$ and 
$\gamma = 2$ or $\gamma = 3$. The other parameters were chosen as in the 
1D test case 
in the next section. The arising patterns are very similar to the
patterns illustrated in the next section. This observation suggests that 
patterns can also evolve when the parameterisation of the autoconversion 
process is nonlinear. However, for an  accretion term
of the form $A_2=a_2q_cq_r$ one can show that even for exponents $\gamma>1$ 
these schemes do not allow Turing instabilities.

\end{enumerate}

\section{Numerical simulations of cloud pattern}
\label{sec:NumericSims}

\emph{Setup:} We carry out 1D and 2D numerical simulations for
investigating the special case ($\beta=2$) of the general cloud model allowing Turing
instabilities, discussed in
section~\ref{sec:CloudSchemePatternFormation}. A pseudo-spectral
method is applied for the numerical solution of the system, see
appendix \ref{sec:Pseudo-spectralMethod}. We choose a domain length of
$L=50$ for the 1D case and a quadratic domain with length of $L=50$
for the 2D case, respectively. In both cases, the domain is
cyclic as assumed in the linear stability analysis above. 
For the 1D simulations, we specify the parameters as follows:\\
$a_1 =1, a_2 = 1, c=5, d=0.1, D_1=10^3$ and $D_2=10^{-1}$.  In this
scenario, the non-trivial stationary  point is asymptotically stable and 
the parabola $p_2$ defined by equation \eqref{eq:quadrPloy}
is negative for wave numbers $q^2= \frac{4 \pi^2}{L^2}n^2$ with
$n \in \{2, \dots, 7\}$. Therefore, these modes give rise for linear
instability and thus lead to Turing instability. 
For the 2D simulations we choose $d=0.13, D_1=10^2, D_2=2.5 \cdot 10^{-2}$. 
The other parameters are like in the 1D case. This slight modification does not 
change the qualitative behaviour. 
As initial condition for both cases
we prescribe the equilibrium of the ODE system with spatial, normally distributed
perturbations with amplitude of order $0.01$. In a first
step, the system~\eqref{eq:CloudScheme_WithPatterns} is simulated,
i.e. there is no rain flux from above ($B=0$). In a second step we
will discuss the impact of the rain flux on the pattern formation in
the simulations. 

\smallskip

\begin{figure}
\centering
\begin{subfigure}[c]{0.49\textwidth}
\includegraphics[width=\textwidth]{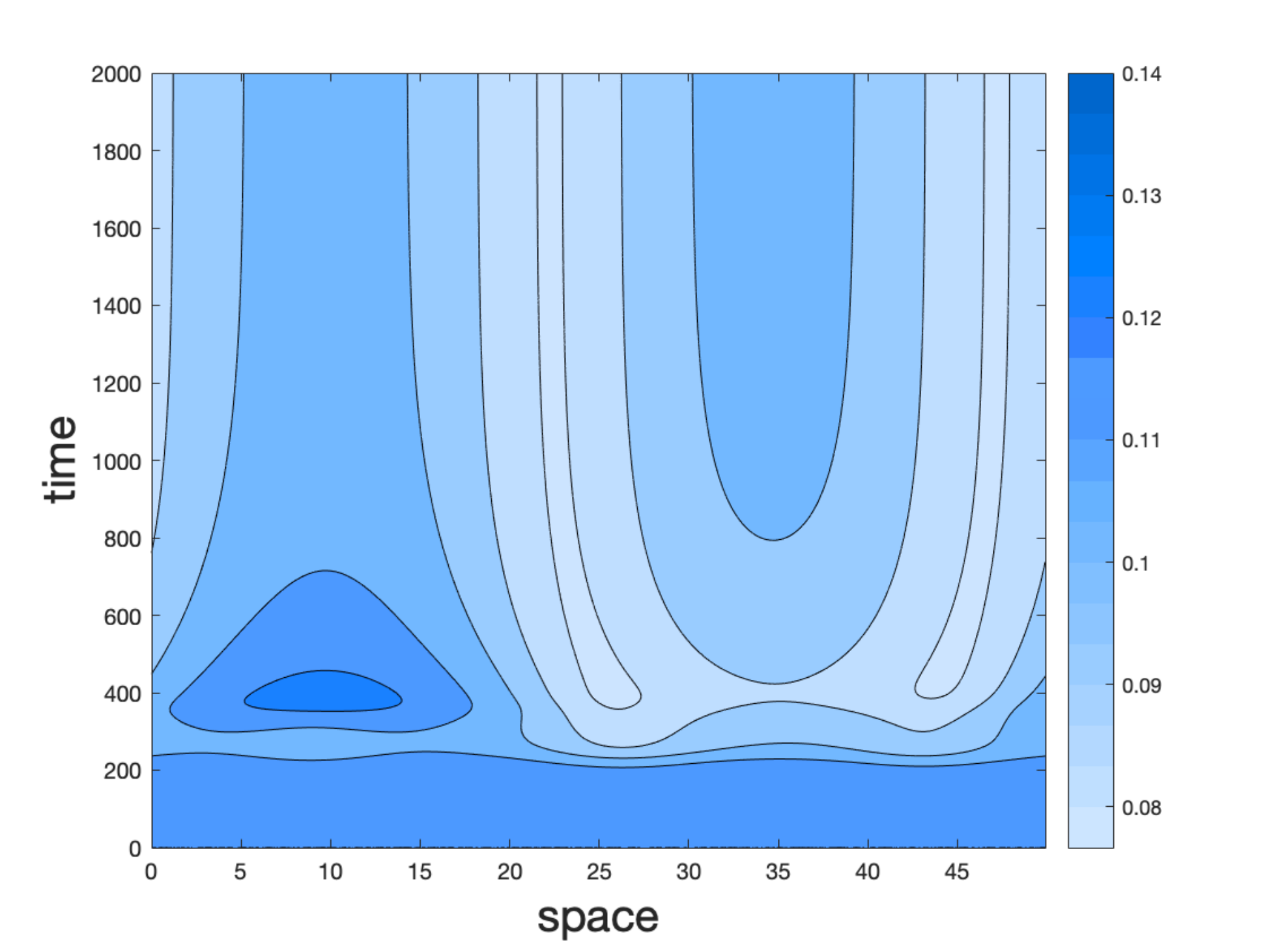}
\subcaption{$q_c$}
\end{subfigure}
\begin{subfigure}[c]{0.49\textwidth}
\includegraphics[width=\textwidth]{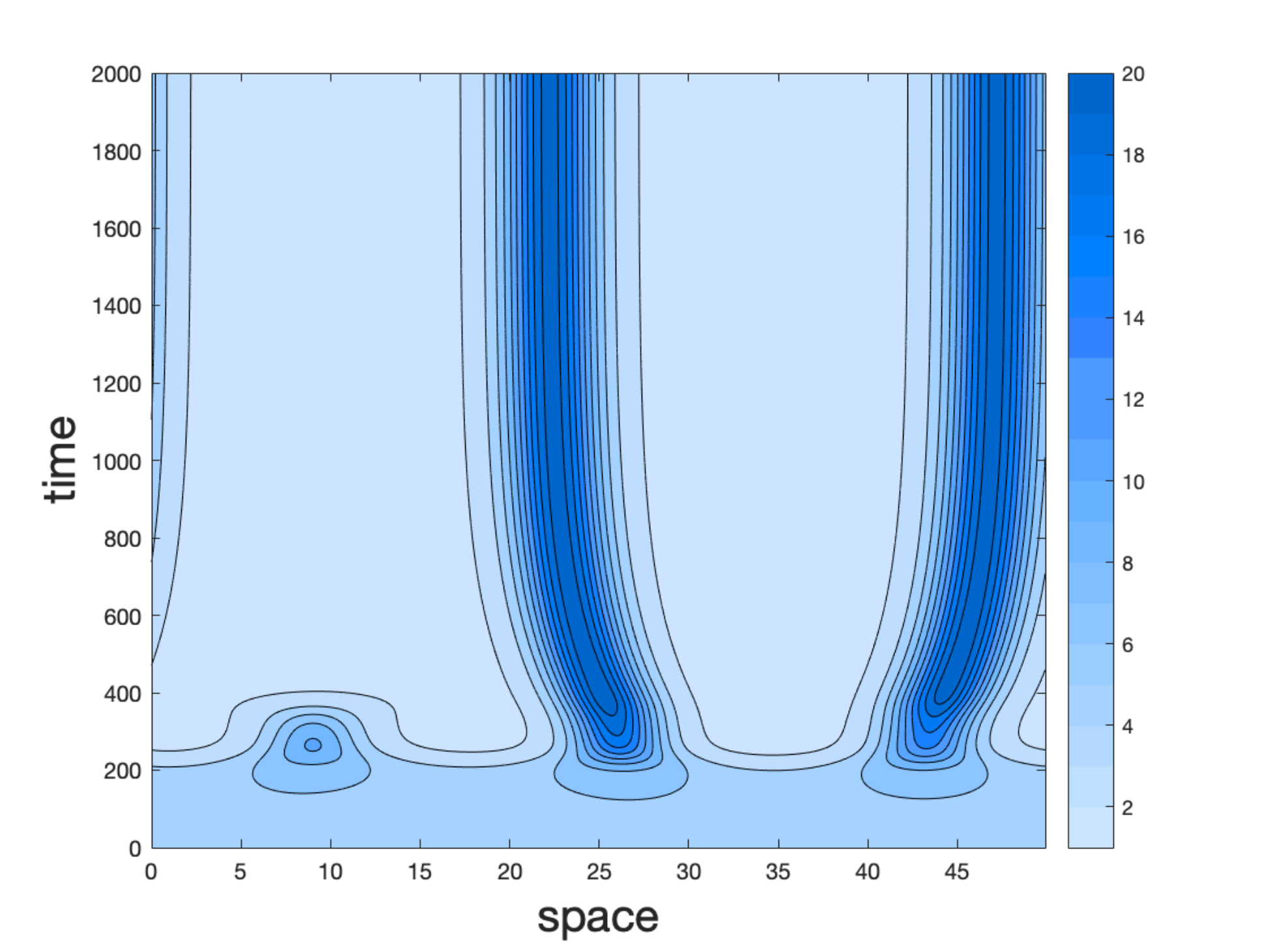}
\subcaption{$q_r$}
\end{subfigure}
\caption{Time evolution of the variables $q_c$ (left panel a) and
  $q_r$ (right panel b) in 1D. The spatial dimension is displayed on
  the  horizontal axis (cyclic domain of length $L=50$).The time is displayed on
  the  vertical axis. At $t\sim 200$ spatial structures form, which finally
  lead to a kind of wavy pattern at the end of the simulation}
\label{fig:pattern1D_qcqr}
\end{figure}

\noindent\emph{Results of 1D simulations without rain flux ($B=0$):}
First, we investigate the numerical simulations in one spatial
dimension. In figure~\ref{fig:pattern1D_qcqr}
the time evolution of the
two variables $q_c$ (left panel) and $q_r$ (right panel) is
shown. The  horizontal axis represents the spatial extension of the 1D domain
(with cyclic boundary conditions), the  vertical axis represents time. 
The values of the cloud variables are represented by the colour code. 
Note, that we always consider dimensionless variables $q_c,q_r$, 
thus the absolute values of these variables have no specific physical meaning.

\begin{figure}
  \centering
  \begin{subfigure}[c]{0.325\linewidth}
    \includegraphics[width=\textwidth]{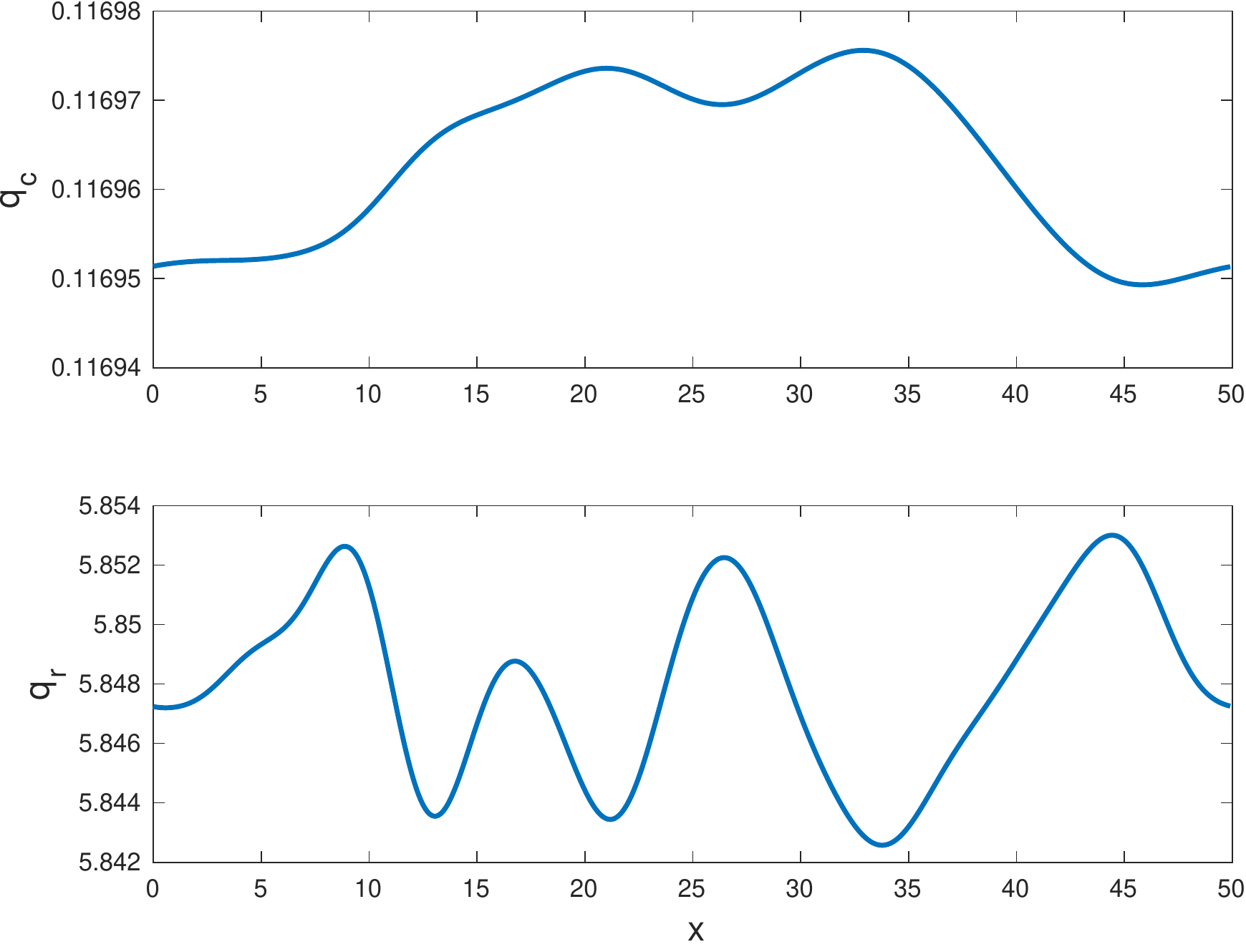}
    \subcaption{$t=20$}
  \end{subfigure}
  \begin{subfigure}[c]{0.32\linewidth}
    \includegraphics[width=\textwidth]{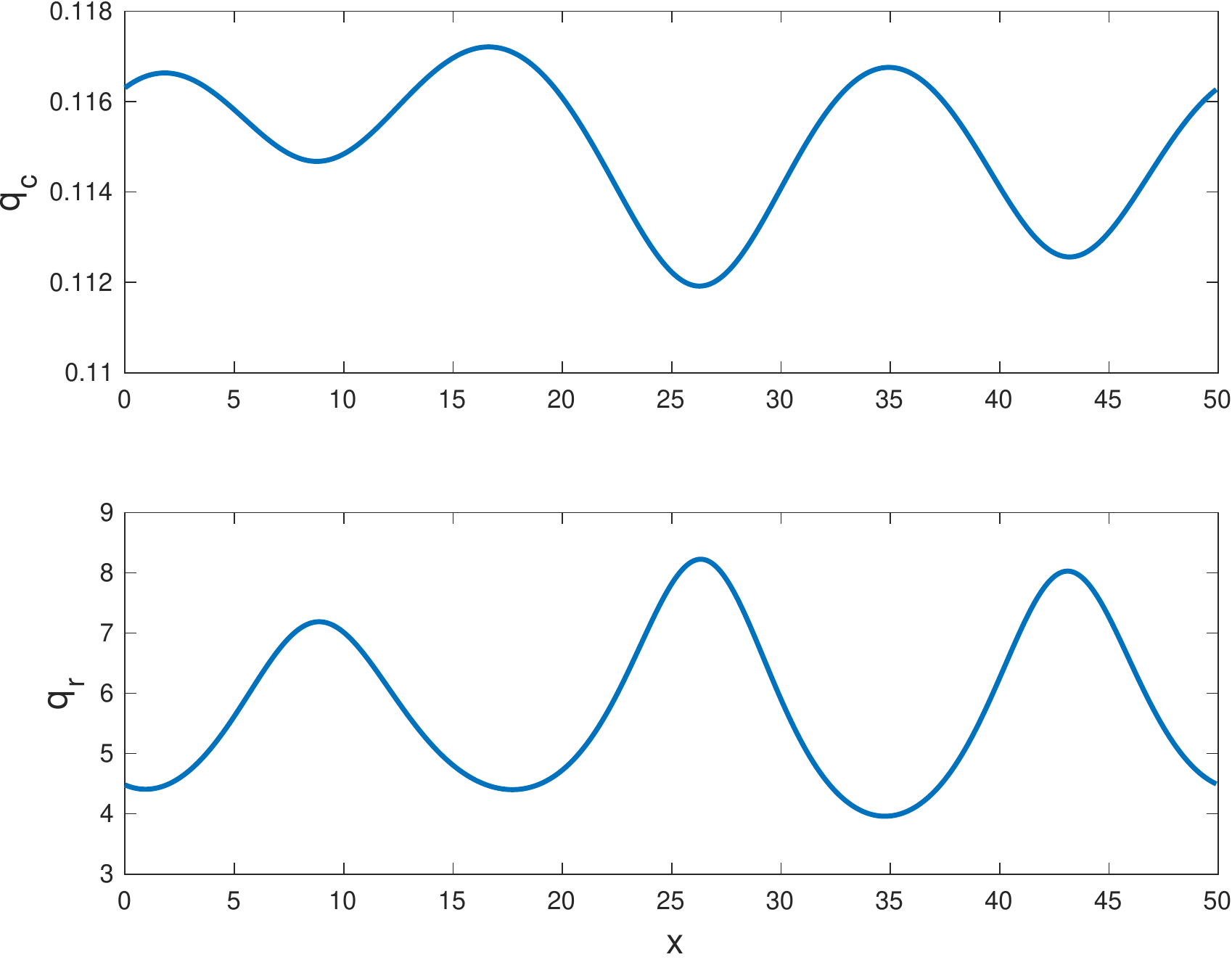}
    \subcaption{$t=200$}
  \end{subfigure}
  \begin{subfigure}[c]{0.32\linewidth}
    \includegraphics[width=\textwidth]{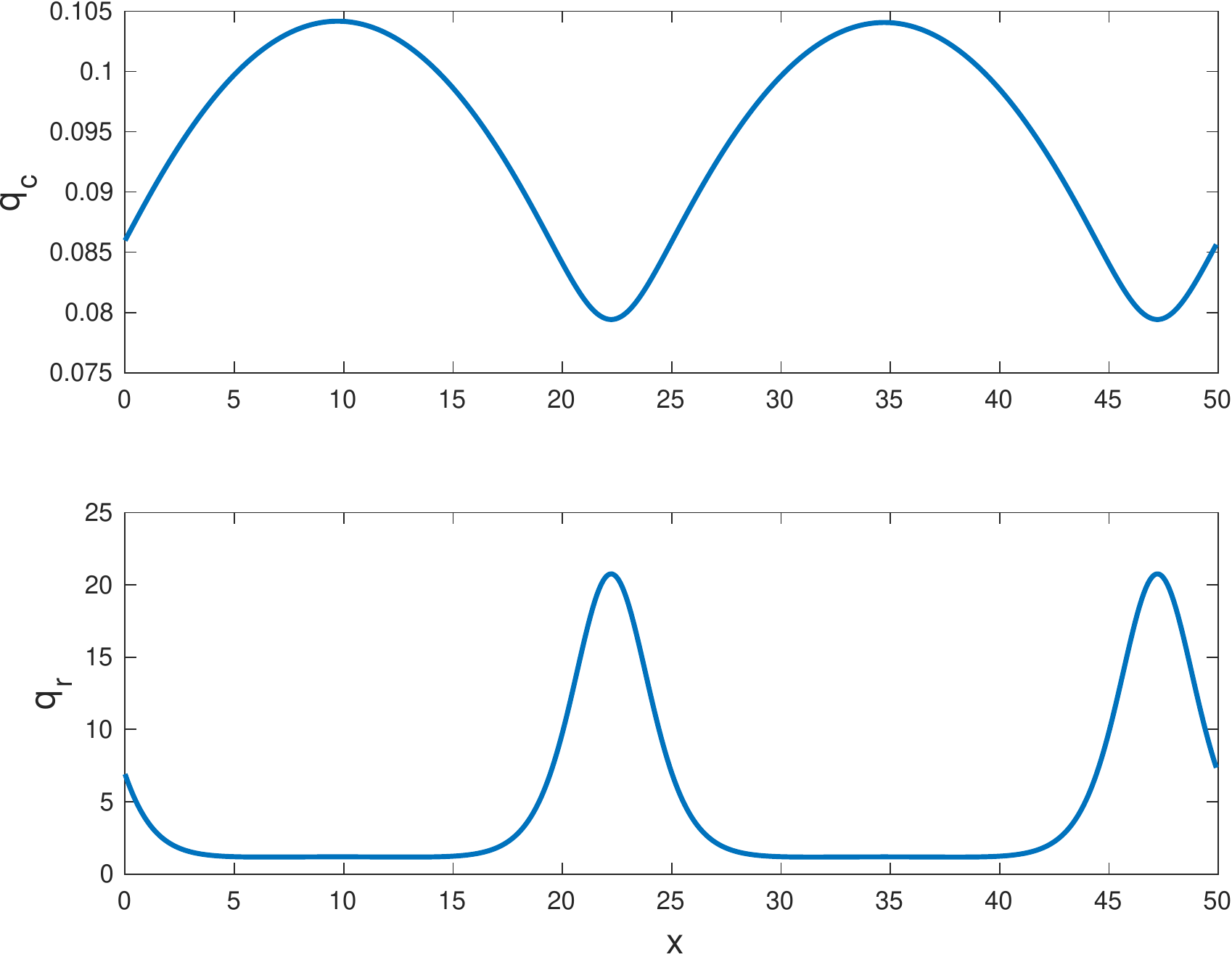}
    \subcaption{$t=2000$}
  \end{subfigure}
  \caption{Spatial variation of the variables $q_c$ (top row) and
    $q_r$ (bottom row) for times $t=20$ (left), $t=200$ (middle), and
    $t=2000$ (right), respectively. Note the different scaling of the
     vertical axes. Actually, at $t=20$ there is almost no variation of
    $q_c,q_r$ visible, whereas at $t=2000$ the change in $q_c,q_r$ is
    obvious.}
  \label{fig:timeslices_qcqr_1D}
\end{figure}

The time evolution clearly shows the formation of spatial structures at times
$t>200$ (in dimensionless time). 
The spatial structure is forming out of the noise, i.e. the
destabilised modes suddenly grow to larger sizes until they are
saturated (and thus stopped) by the nonlinear terms. Their spatial
distribution slightly changes during time until at around $t\sim 1200$ the
situation is consolidated, i.e. the pattern stays quite stationary. In
figure~\ref{fig:timeslices_qcqr_1D} the simulations at times
$t=20/200/2000$ are shown. Here, the evolution can be seen clearly as
well as the final ``wavy'' structure at $t=2000$. Note, that the
variables $q_c$ and $q_r$ have contrary behaviour: for high values of
$q_c$ the rain variable $q_r$ is quite small and vice versa. This can
be explained by the collision terms, which act as in a generalised
predator-prey system. If the predator population (i.e. the rain) is
small, the cloud water survives and grows to larger values due to
condensation only, since autoconversion is weak. 
If the rain becomes larger, it reduces the prey
(the cloud water) due to collisions.

Using Fourier analysis (not shown) we see that only a part of the
Fourier spectrum has reasonable amplitudes, whereas higher modes are of
very low amplitudes. However, we do not see the distinct spectrum as
predicted by linear stability theory. The reason for this is the
nonlinear interaction of the different modes, which 
leads to non-vanishing amplitudes of modes which are stable according to the linear stability analysis. 
Nevertheless, we see that only a small part of the
Fourier spectrum is present in the simulations.

\smallskip
\noindent\emph{Results of 2D simulations without rain flux ($B=0$):}
In a second simulation we use the 2D setup with white noise as before to
investigate pattern formation in a 2D domain.

Qualitatively, we see the same behaviour for the 2D simulations of a
quadratic domain of length $L=50$ with cyclic boundary
conditions. After a short time, the simulation leads to growing
unstable modes, which are then saturated by nonlinear terms in the
model; these modes form spatial structures, which change only slightly
over time until they stay stationary. Thus, pattern formation due to
Turing instabilities can be observed as expected from theory.
The structures in cloud water $q_c$ are less pronounced than in the rain water $q_r$. 
Nevertheless, the spatial patterns remain stationary, even for longer times.

\begin{figure}
\centering
\begin{subfigure}[c]{0.49\textwidth}
\includegraphics[width=\textwidth]{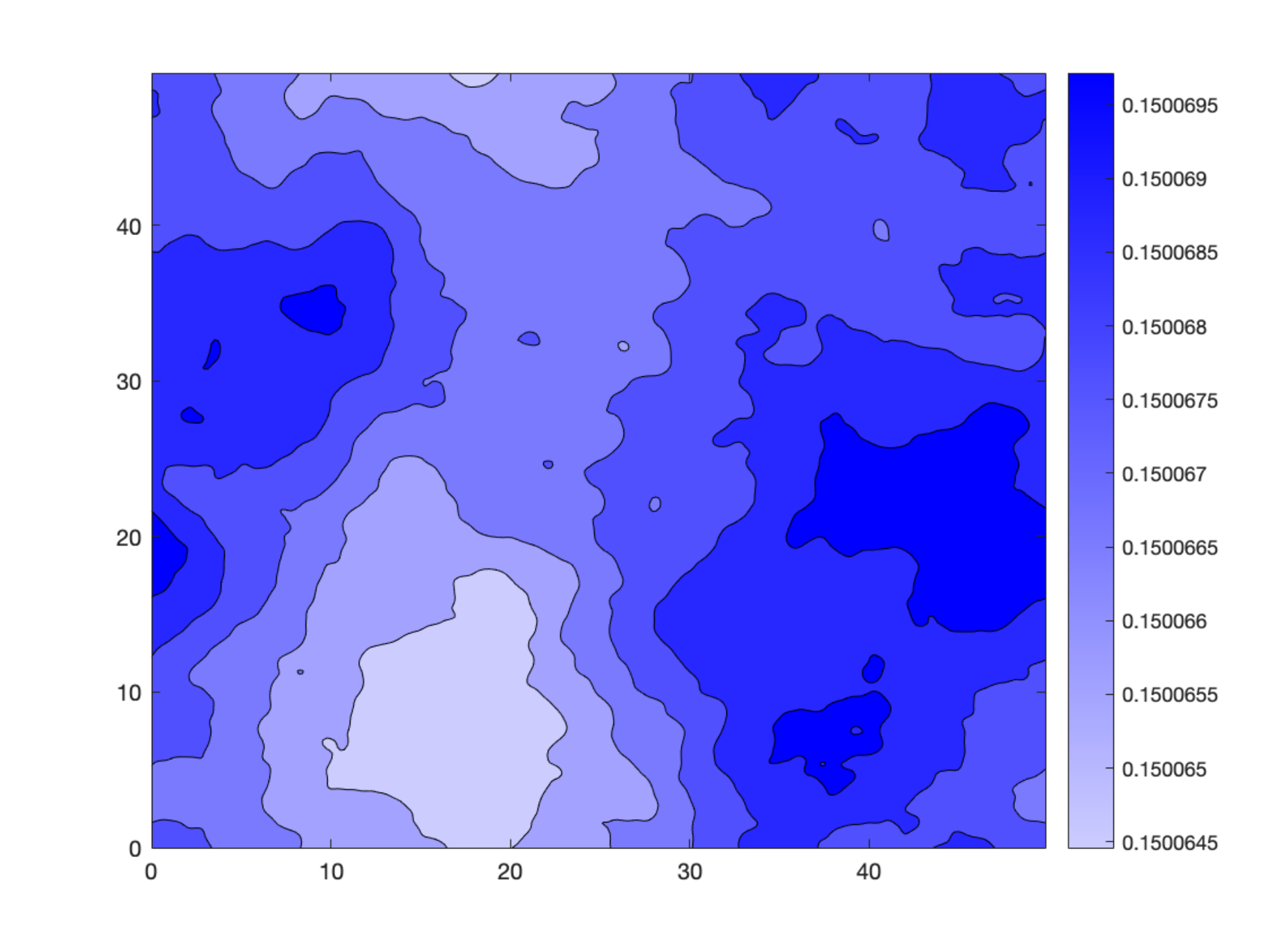}
\subcaption{$t=1$}
\end{subfigure}
\begin{subfigure}[c]{0.49\textwidth}
\includegraphics[width=\textwidth]{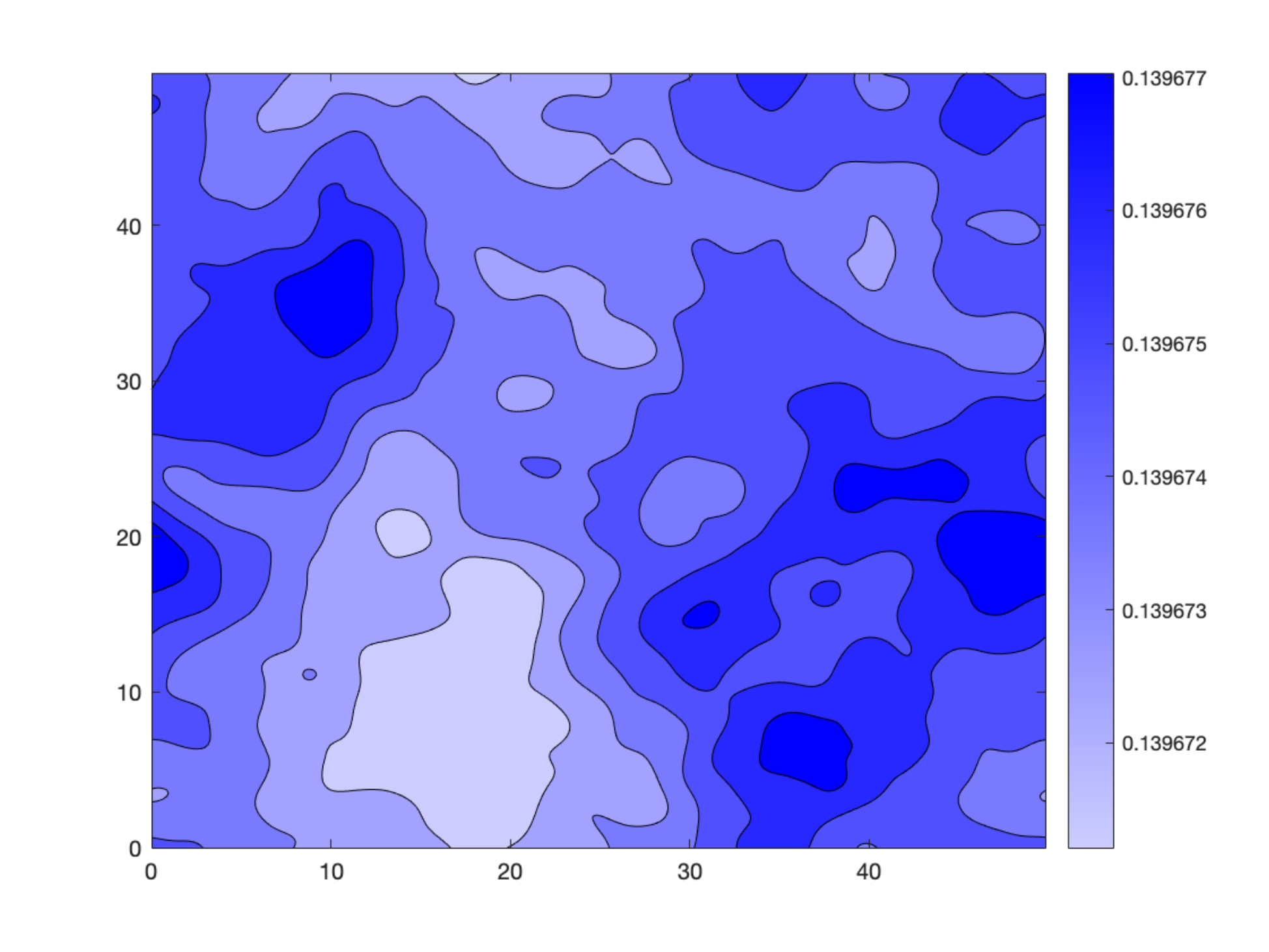}
\subcaption{$t=10$}
\end{subfigure}
\begin{subfigure}[c]{0.49\textwidth}
\includegraphics[width=\textwidth]{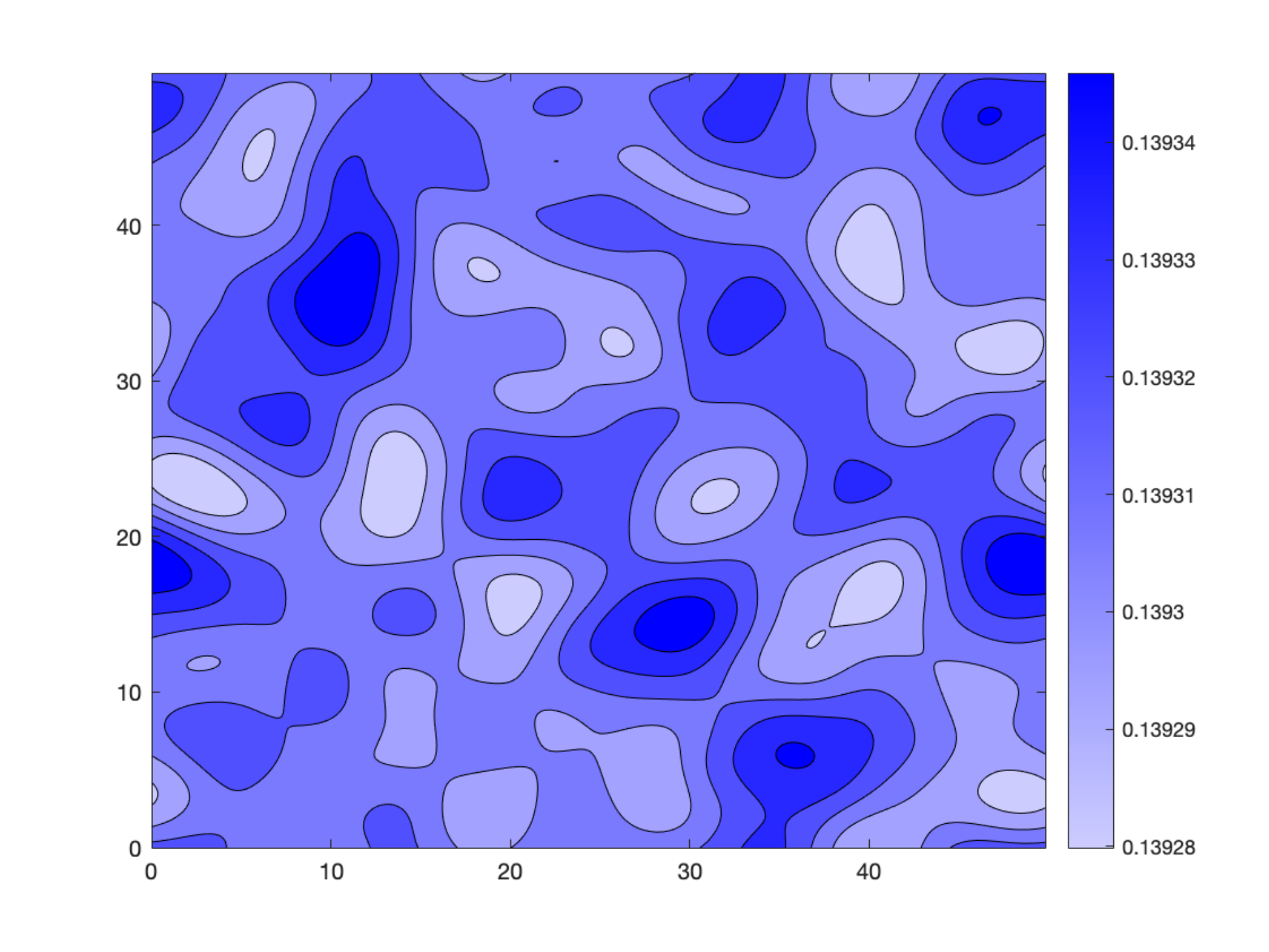}
\subcaption{$t=60$}
\end{subfigure}
\begin{subfigure}[c]{0.49\textwidth}
\includegraphics[width=\textwidth]{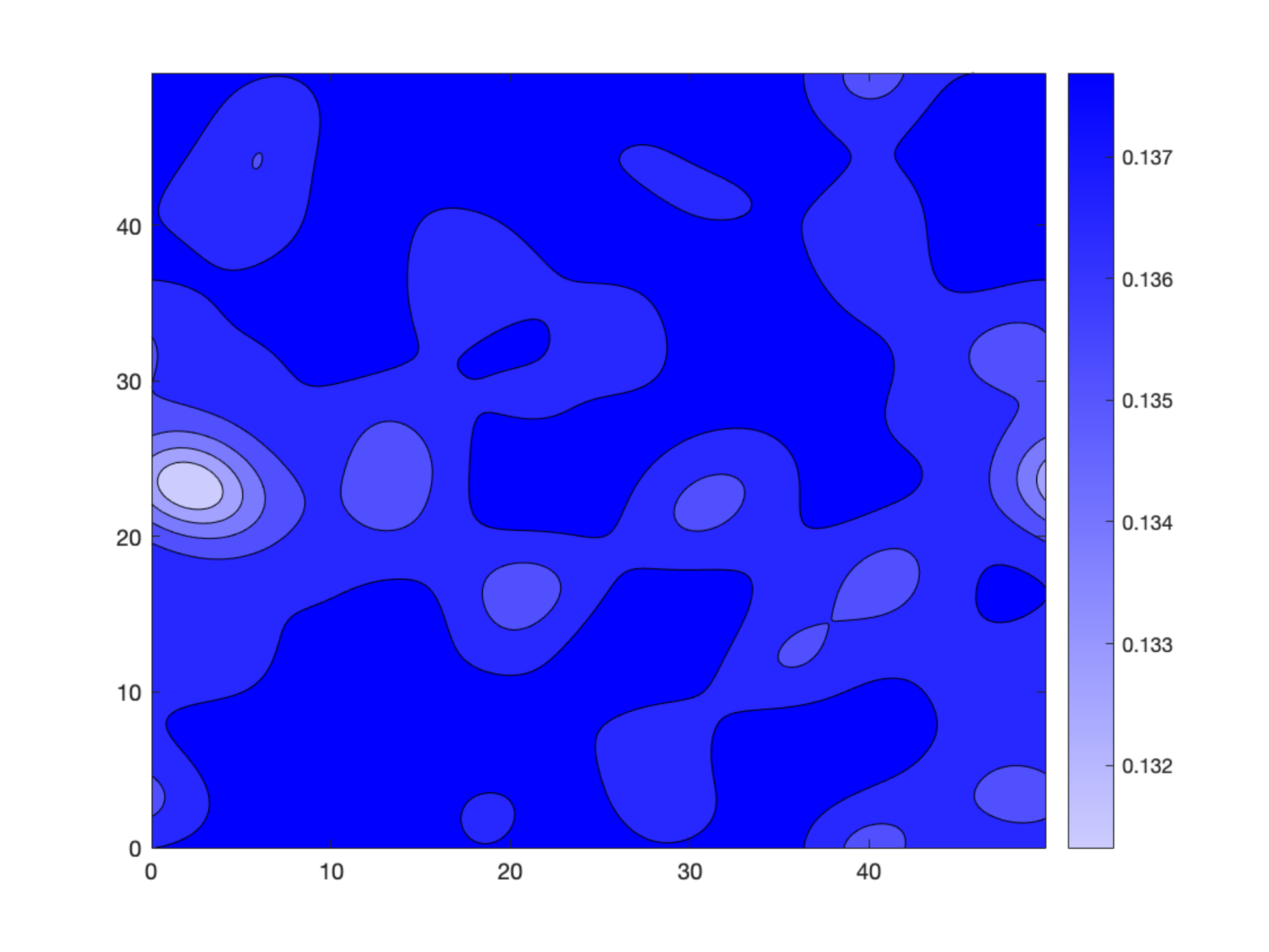}
\subcaption{$t=120$}
\end{subfigure}
\caption{Spatial distribution of cloud water $q_c$ in 2D for different
  simulation times ($t=1/10/60/120$). Note, that the pattern is
  already forming at times $t\sim 10$. For longer times, the pattern
  stays stationary although the absolute variation of the cloud
  water variable is very small over the whole 2D domain.}
\label{fig:pattern2D_qc}
\end{figure}

\begin{figure}
\centering
\begin{subfigure}[c]{0.49\textwidth}
\includegraphics[width=\textwidth]{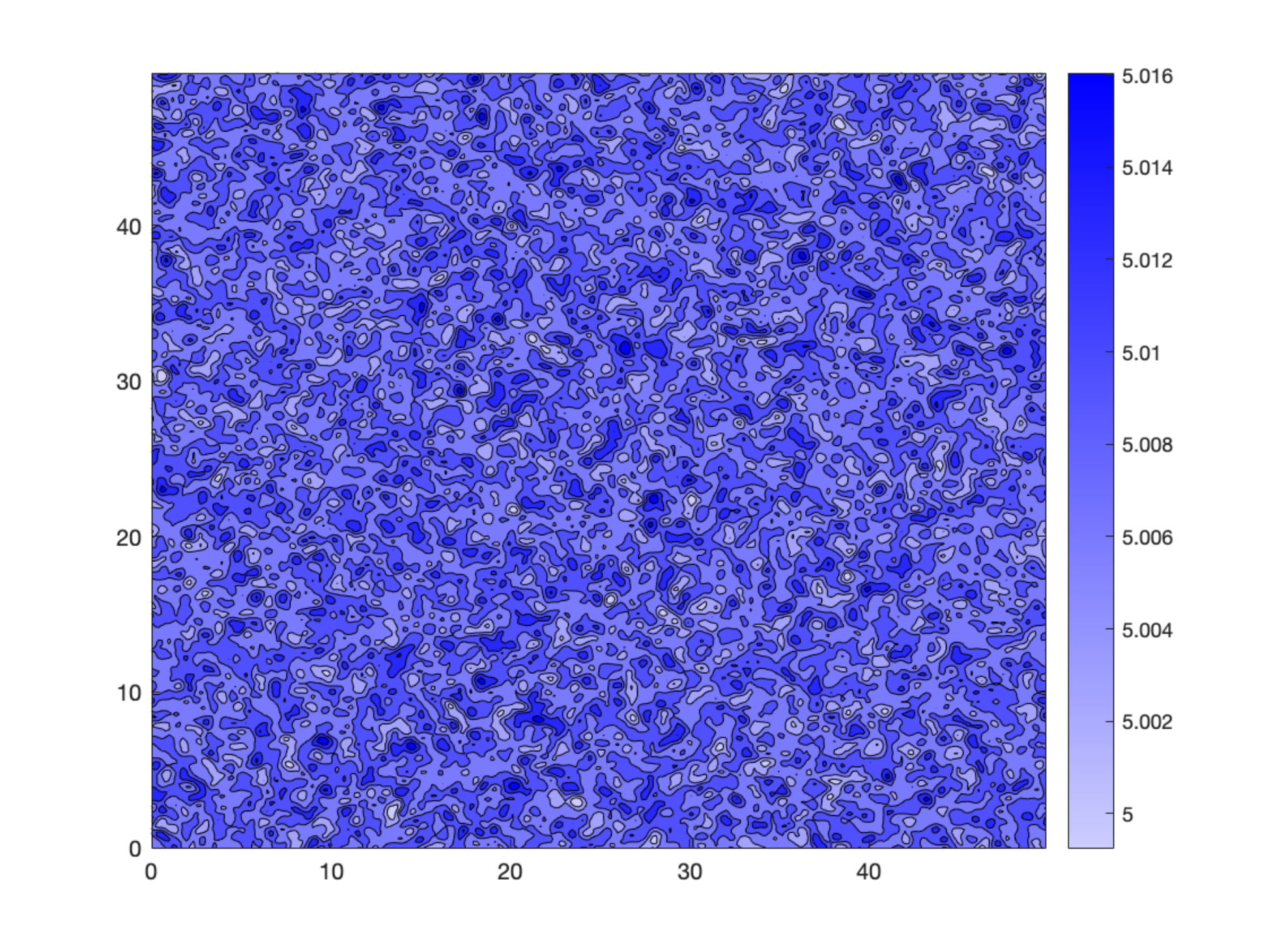}
\subcaption{$t=1$}
\end{subfigure}
\begin{subfigure}[c]{0.49\textwidth}
\includegraphics[width=\textwidth]{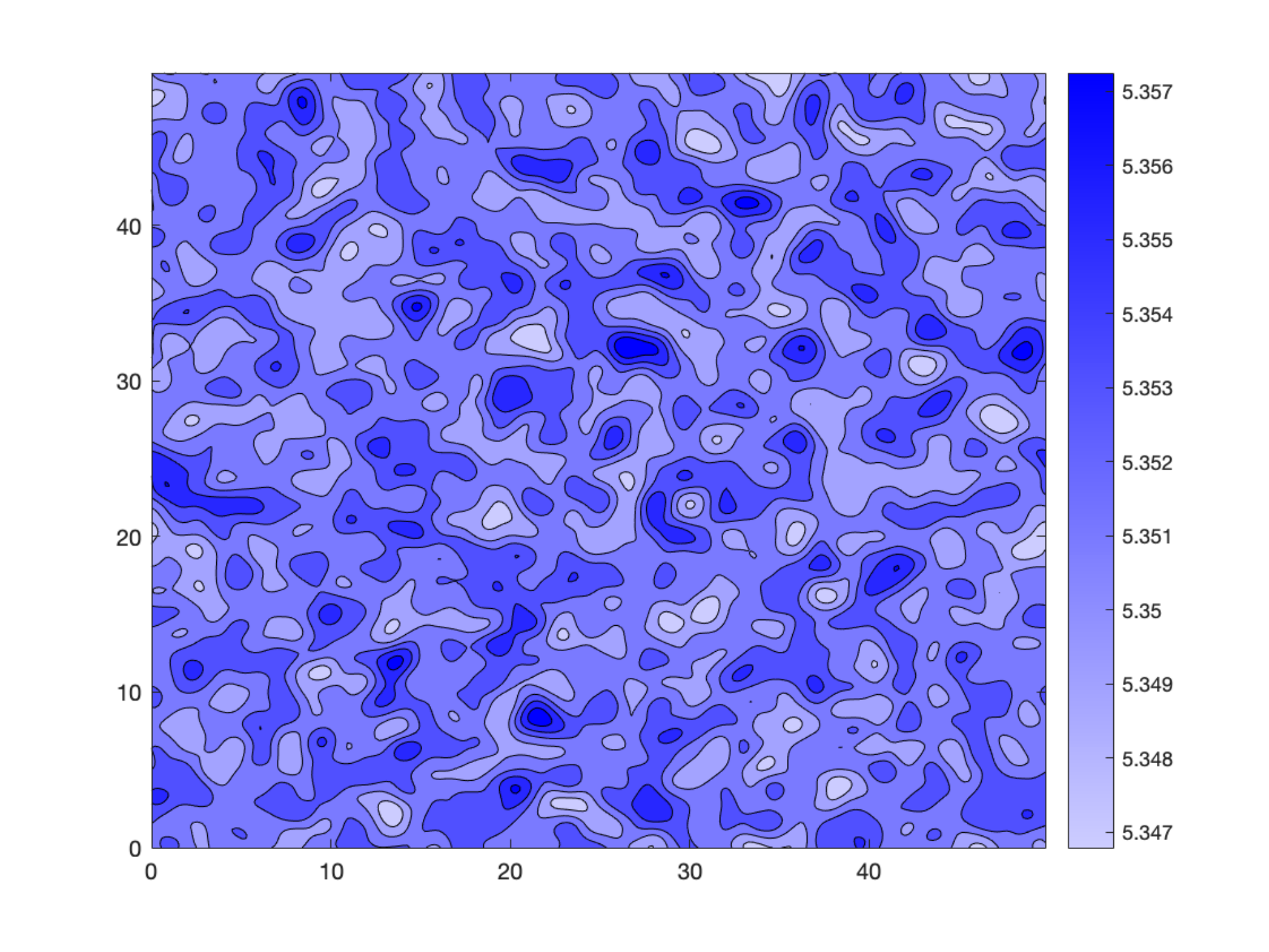}
\subcaption{$t=10$}
\end{subfigure}
\begin{subfigure}[c]{0.49\textwidth}
\includegraphics[width=\textwidth]{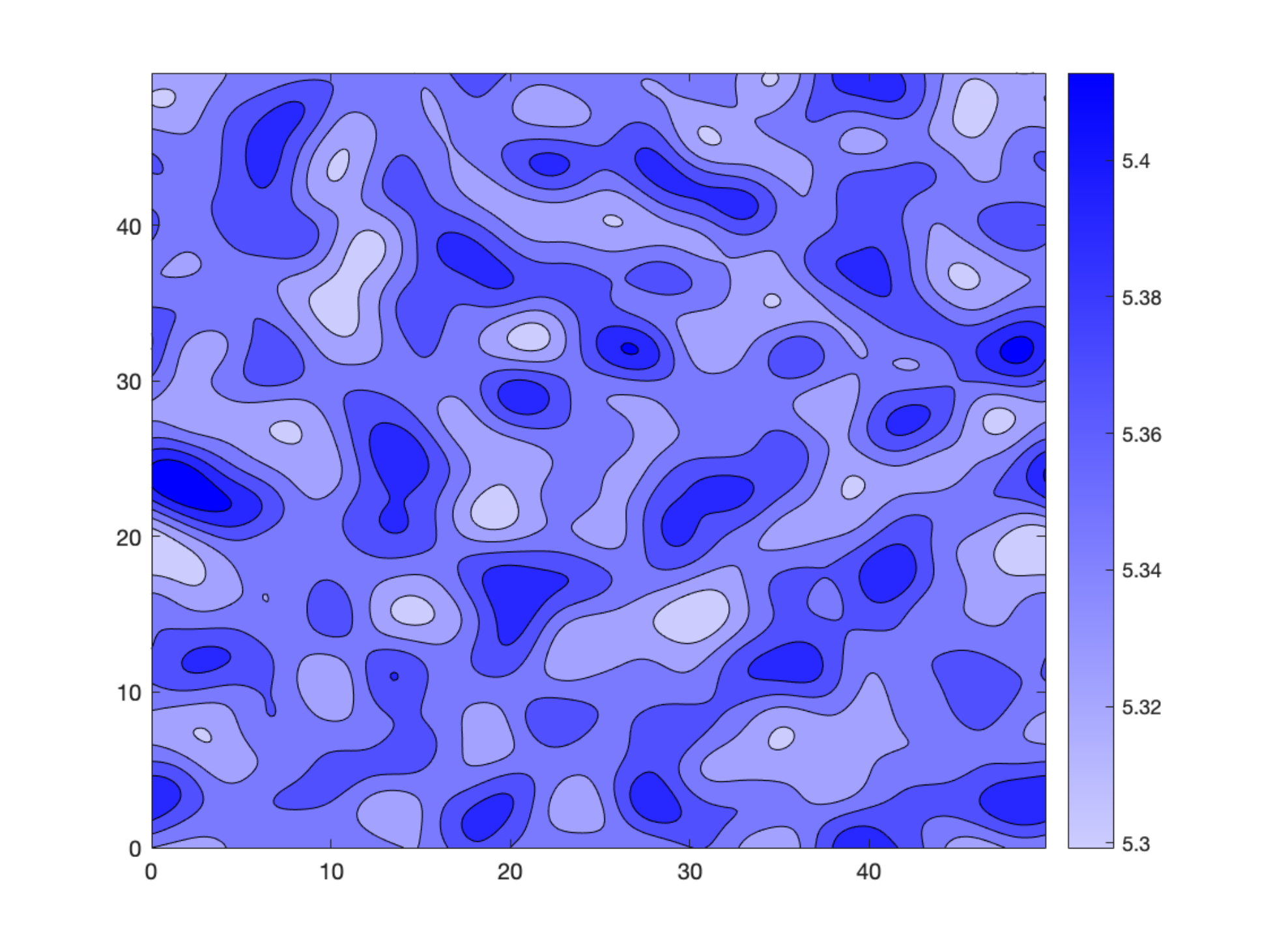}
\subcaption{$t=60$}
\end{subfigure}
\begin{subfigure}[c]{0.49\textwidth}
\includegraphics[width=\textwidth]{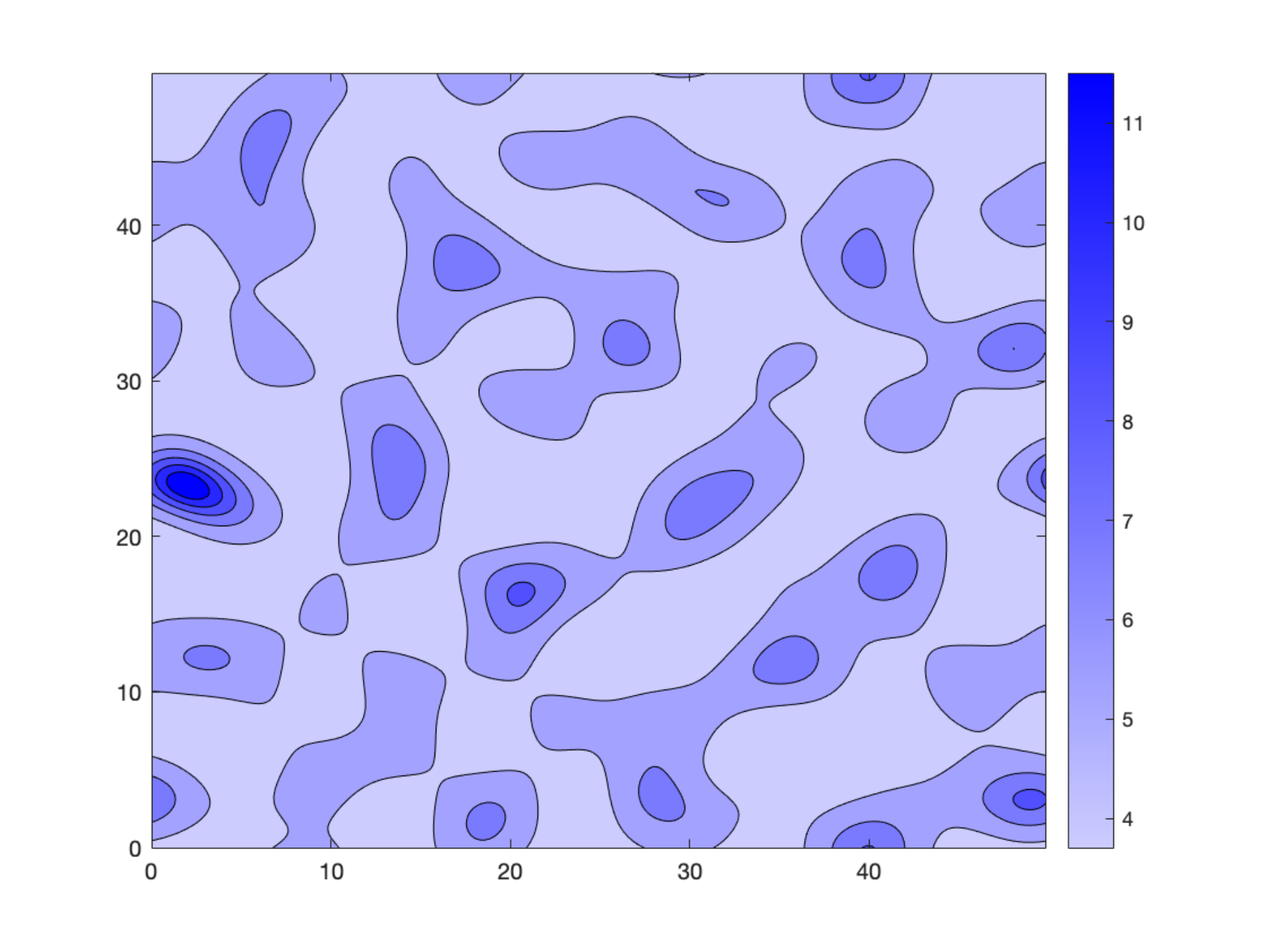}
\subcaption{$t=120$}
\end{subfigure}
\caption{Spatial distribution of rain water $q_r$ in 2D for different
  simulation times ($t=1/10/60/120$). The spatial structure for this
  variable is more pronounced than for the cloud water, i.e. the
  spatial variation of $q_r$ is quite large.}
\label{fig:pattern2D_qr}
\end{figure}

\smallskip

\noindent\emph{Results of 1D simulations for including the rain flux ($B>0$):}
In a last series of simulations we investigate the impact of the rain
flux $B$, which was set to zero in section~\ref{sec:CloudSchemePatternFormation}
for simplification of the
analysis. Using the fixed parameters $c,a_1,a_2,d$ we can calculate
the roots of the cubic polynomial determining the stationary points. One real
root can be calculated as $q_{re}=u+v+\frac{B}{3d}$, 
see \eqref{sec6:real_root}. The bifurcation
value can be calculated as $B_1\sim1.10521$. Actually, we additionally
find that the eigenvalues $\sigma_i$ (for $i=1,2$)  always have negative 
real part for
$0<B<B_1$.
Thus, the nontrivial
stationary point is always asymptotically stable for rain flux from
above in the relevant parameter range.

As described in section~\ref{sec:CloudSchemePatternFormation}, the
entries $a_{11}<0,a_{21}>0$ of the Jacobi matrix are constant, the
entries $a_{12}(B),a_{22}(B)$ depend on the rain flux. For fulfilling
the criterion for Turing instabilities~\eqref{eq:unstableDiffusionSystem}, 
the sign of entry $a_{22}$ decides about the existence
or non-existence of instabilities. For values
$B<B_2\sim 0.137$ we obtain
$a_{22}(B)>0$ (i.e. Turing instability is possible), whereas for $B>B_2$
the entry is negative.

\begin{figure}
\centering
\begin{subfigure}[c]{0.49\textwidth}
\includegraphics[width=\textwidth]{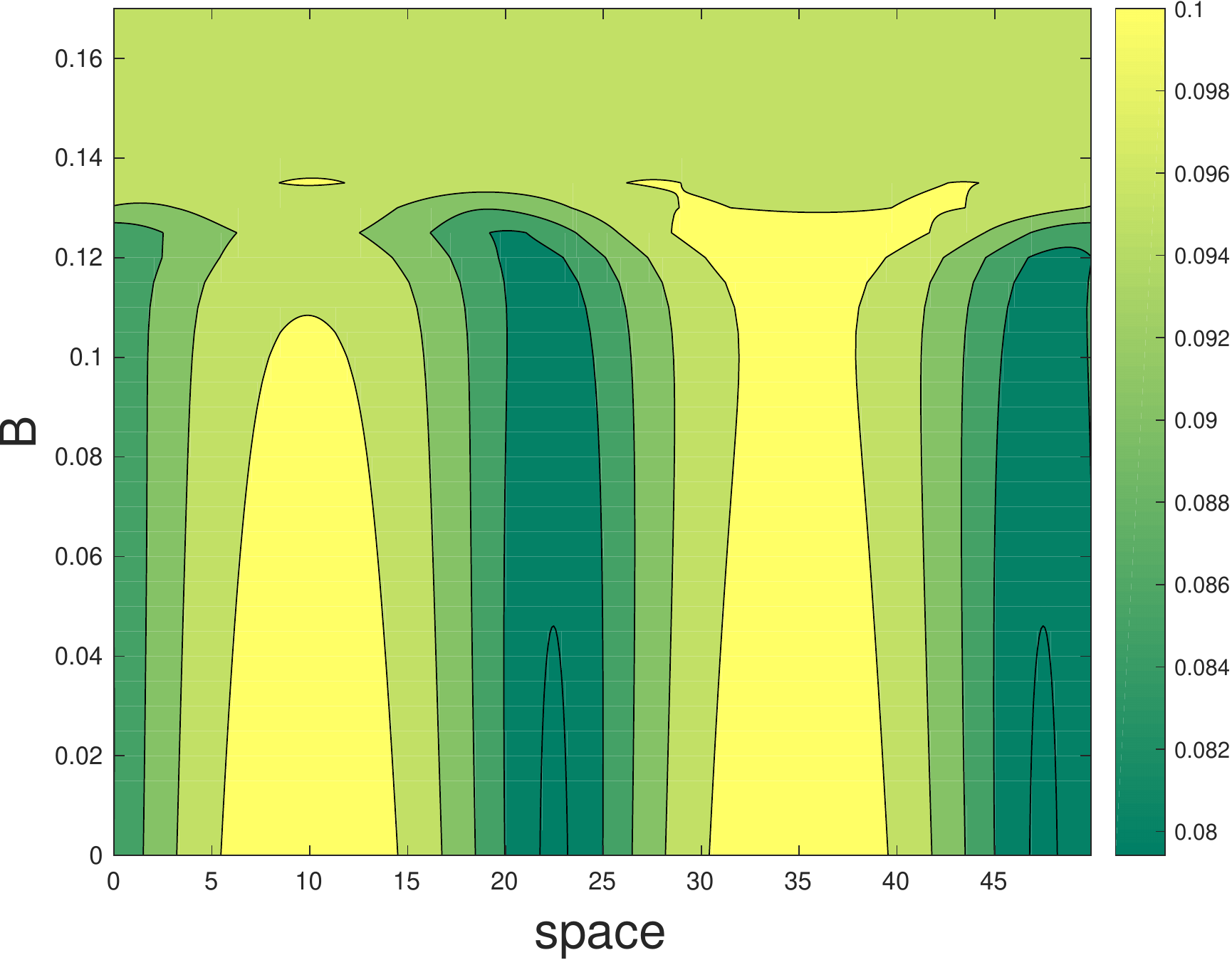}
\subcaption{$q_c$}
\end{subfigure}
\begin{subfigure}[c]{0.49\textwidth}
\includegraphics[width=\textwidth]{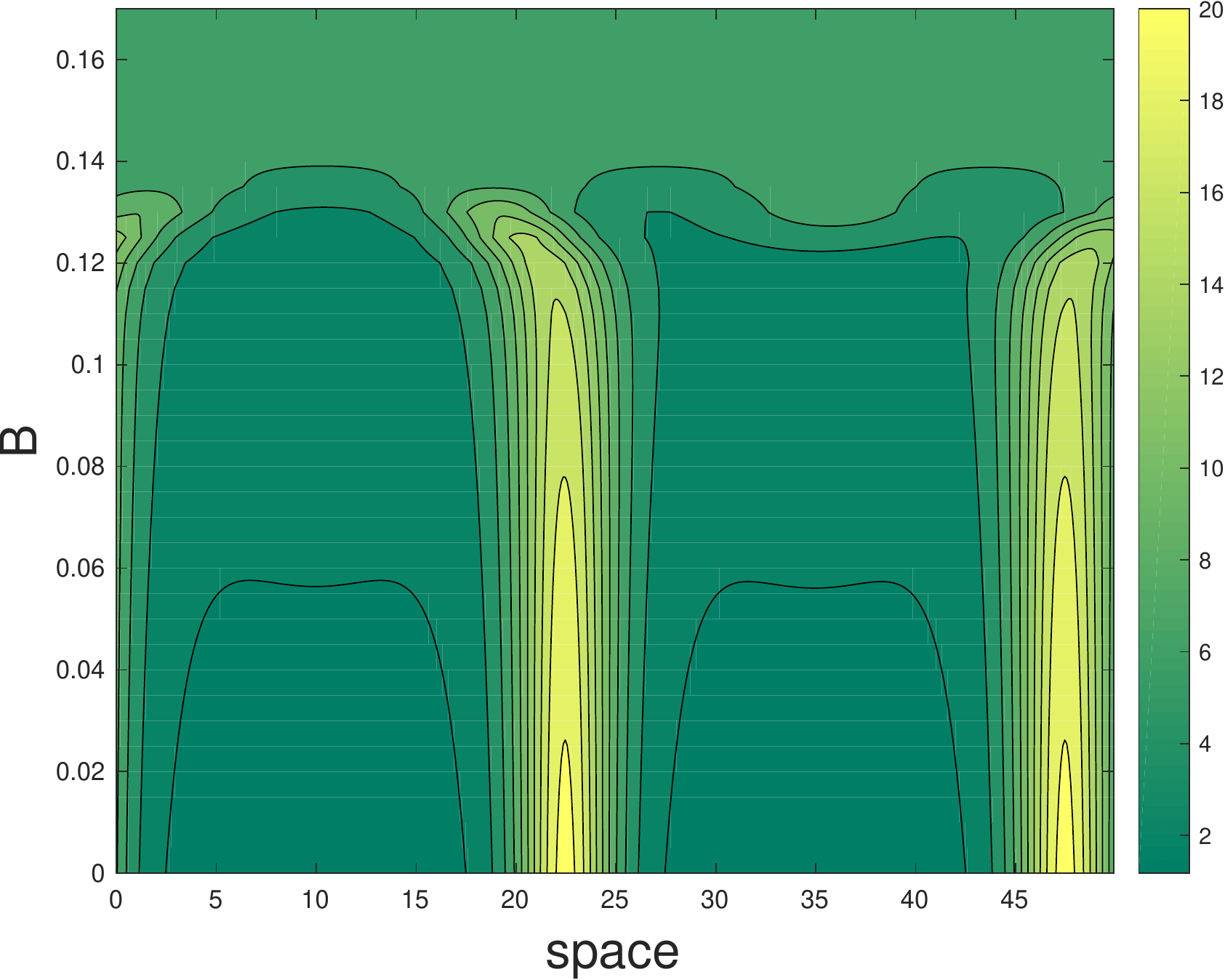}
\subcaption{$q_r$}
\end{subfigure}
\caption{The cloud variables $q_c$ (left panel a) and $q_r$ (right panel b) 
after 2000 time units for different values of the rain flux $B$ in 1D simulations. 
The x-axis displays the spatial direction, the y-axis represent different values 
of the rain flux $B$. The values of the cloud variables are given by the colour code.}
\label{fig:Turing1D_B_qcqr}
\end{figure}

We confirm these findings with a series of numerical simulations
using different values $0<B<0.17$ for the set of parameters as specified 
at the beginning of the section. 
As predicted, for values $B<0.137$ we 
find Turing instabilities, whereas for $B>0.137$ there are no Turing
instabilities. In Figure~\ref{fig:Turing1D_B_qcqr} the simulations at
time $t=2000$ (i.e. steady state) depending on the parameter $B$ are
shown. The absolute values of the pattern in $q_c$ and $q_r$
slightly vary with changing $B$; however, the quality of the pattern remains
the same until values $B\sim B_2$ are reached. Passing this values, a
homogeneous state in both variables can be seen and no pattern formation occurs.
Note, that the boundary $B_2$ is not sharp, since the occurrence of 
the Turing instability depends also on the values of $D_1,D_2$. Due to the 
large ratio of these coefficients, the transition in the simulations 
is very close to $B_2$.

\section{Summary and conclusion}
\label{sec:SummaryConclusion}

In this study we investigate a generic  cloud model for
pure liquid clouds on the possibility of Turing instabilities for
forming spatial patterns. This kind of investigation is carried out for 
the first time for 
cloud pattern formation. For the theoretical and numerical investigations, 
the generic cloud model 
formulated in the former study \cite{rosemeier_etal2018} is extended
by diffusion terms, consisting of Laplacians in spatial
directions. The model is analysed using stability theory for the 
linearisation around
the steady states of the underlying two-dimensional ODE system. 
Analytical conditions for the 
existence of Turing instabilities can be determined.
Since the model contains quite complex nonlinear terms
with several parameters determining the overall quality of the steady states,
it is very hard to find general conditions for the existence of Turing instabilities.
However, the generic model always admits a trivial stationary point in terms of ``no clouds, 
just rain falling through the layer''. This stationary point could be either stable or unstable,
depending on the set of parameters. However, even in the stable case, 
the steady state cannot be destabilised by diffusion; thus, this state does not 
admit Turing instabilities. 
In addition, we can specify a class of models, which do not allow Turing instabilities
at all. This class is characterised by a linear autoconversion term $A_1=a_1q_c$ and also 
a linear contribution of cloud water $q_c$ in the accretion term $A_2$. Well-known 
cloud models such as the standard COSMO cloud scheme
(\cite{cosmo_doc_physical_parameterization})  
or the research model by \cite{wacker1992} belong to this class. 
However, we can also provide a general class of cloud models, which allow Turing
instabilities. If the exponents in the accretion parameterisation are chosen to be larger 
than 1 (i.e. $\beta=\beta_c=\beta_r>1$), the model can allow Turing instabilities 
and thus pattern formation. These theoretical findings could be confirmed by numerical 
simulations in one and two spatial dimensions. The inclusion of rain flux from above 
turns out to be an additional restriction for the instabilities. 
This is investigated for the special case $\beta=2$.
If the rain flux becomes too
large (i.e. if it surpasses a certain threshold $B>B_2$), the criterion for the existence of 
Turing instabilities is violated, as can be seen also in a series of 1D numerical simulations.
This observation leads to the interpretation  that collision processes in combination 
with the sedimentation of cloud particles play the major role for pattern formation; 
only if these processes can interact in a proper nonlinear way, Turing instabilities 
are possible. A strong rain flux from above can prevent the formation of cloud patterns;
this might be explained by stronger collision terms, which finally almost extinct the cloud
droplet population, thus diffusion can not counteract this process.

We can conclude that the generic cloud model admits Turing instabilities in special cases.
However, several standard cloud models  used in research and operational weather forecasts
do not admit pattern formation due to Turing instabilities. It is still unknown how patterns 
in clouds form, especially, which processes lead to the emergence of cloud structures.
The use of diffusion terms is motivated by the parameterisation of subgrid scale 
processes; this approach might be too simplistic for representing the underlying processes
in a meaningful way. On the other hand, it is conceivable that pattern formation 
in clouds
is dominated by Turing instabilities; in this case it might be a major drawback to use 
cloud models, which do not allow this type of pattern formation. Since pattern formation
in clouds is far away from being understood, this observation has to be taken into account, 
and it might have an impact on the choice of cloud models for further investigations.
For the investigation of cloud patterns, more theoretical studies are needed for a better
understanding of the underlying processes and their interaction, which leads to the 
emergence of cloud structures.

\smallskip

\noindent\textbf{Acknowledgement:} We thank Maria Lukacova and Manuel Baumgartner 
for fruitful discussions. 
We acknowledges support of the Transregional Collaborative Research Center
SFB/TRR 165 ``Waves to Weather'', funded by the ``Deutsche Forschungsgemeinschaft'' (DFG), 
within the sub-project ``Structure Formation on Cloud Scale and Impact on Larger Scales'' 
(Project A2).

\begin{appendix}
  
  \section{Jacobian of ODE system \eqref{eq:generic_generic_model}}

The Jacobian $Df$ of the generic cloud model (without diffusion) can be calculated as

\begin{equation}
Df =
      \begin{pmatrix} 
  c-\gamma a_1 q_c^{\gamma-1}-a_2\beta_c q_c^{\beta_c-1}q_r^{\beta_r}   &  
  -a_2\beta_r q_c^{\beta_c}q_r^{\beta_r-1}\\
  a_1\gamma q_c^{\gamma-1} +
  a_2\beta_c q_c^{\beta_c-1}q_r^{\beta_r}& 
  a_2\beta_r q_c^{\beta_c}q_r^{\beta_r-1} -d\zeta q_r^{\zeta-1}  \\  
  \end{pmatrix}. 
\end{equation}

In case of the trivial equilibrium state 
$q_{ce}=0,~ q_{re}=\left( \frac{B}{d}\right)^{\frac{1}{\zeta}}$
the Jacobian reduces to 

\begin{equation}
\label{eq:JacobianTrivialState}
    \begin{pmatrix}
    a_{11} & 0 \\
    a_{21} & -d\zeta \left( \frac{B}{d}\right)^{\frac{\zeta-1}{\zeta}}
    \end{pmatrix}
\end{equation}
Remember, that the Jacobian of the trivial state is only defined for 
values $\gamma\ge 1$, $\beta_c\ge 1$; otherwise 
the partial derivatives with respect to $q_c$ do not
exist and 
the initial value problem is potentially not uniquely solvable, since the 
right hand side of the ODE system 
\eqref{eq:generic_generic_model}
is not Lipschitz continuous.
For the entries in the matrix \eqref{eq:JacobianTrivialState}, 
we have to discriminate between different cases. First, we determine the 
value of $a_{11}$.
\begin{itemize}
    \item For  $\gamma=1, \beta_c = 1$ we obtain $a_{11}=c-a_1-a_2 \qty(\frac{B}{d})^{\frac{\beta_r}{\zeta}}$. In this case, 
    the trivial stationary point can be stable ($c<a_1 + a_2 \qty(\frac{B}{d})^{\frac{\beta_r}{\zeta}}$) 
    or unstable ($c>a_1 + a_2 \qty(\frac{B}{d})^{\frac{\beta_r}{\zeta}}$).
    \item For  $\gamma=1, \beta_c > 1$ we obtain $a_{11}=c-a_1$. In this case, 
    the trivial stationary point can be stable ($c<a_1$) 
    or unstable ($c>a_1$).
    \item For  $\gamma>1, \beta_c = 1$ we obtain $a_{11}=c-a_2 \qty(\frac{B}{d})^{\frac{\beta_r}{\zeta}}$. In this case, the trivial stationary point can be stable 
    ($c< a_2 \qty(\frac{B}{d})^{\frac{\beta_r}{\zeta}
    }$) or unstable ($c> a_2 \qty(\frac{B}{d})^{\frac{\beta_r}{\zeta}}$). 
    \item For  $\gamma>1, \beta_c > 1$ we obtain $a_{11}=c$. In this case, the trivial stationary point 
    is always unstable.
\end{itemize}
Second, the entry $a_{21}$ is investigated.
\begin{itemize}
    \item For $\gamma=1,\beta_c=1$ we obtain 
    $a_{21}=a_1+a_2\left( \frac{B}{d}\right)^{\frac{\beta_r}{\zeta}}$.
    \item For $\gamma=1,\beta_c>1$ we obtain
    $a_{21}=a_1$.
    \item For $\gamma>1,\beta_c=1$ we obtain
    $a_{21}=a_2\left( \frac{B}{d}\right)^{\frac{\beta_r}{\zeta}}$.
    \item For $\gamma>1,\beta_c>1$ we obtain $a_{21}=0$.
\end{itemize}
In any case, the entry $a_{21}$ does not affect the stability of the trivial stationary point.

  \section{Pseudo-spectral method}
  \label{sec:Pseudo-spectralMethod}

  The pseudo-spectral method is applied to the following type of
  semilinear equations
 \begin{equation}
\label{eq:Pseudospektral-Verfahren-AusgansDgl}
\dv{x}{t} = L \qty(x) + R \qty(x).
\end{equation}
where $L$ is a linear operator and $R$ a nonlinear operator. 
The model equation \eqref{eq:System} represents such a system, the linear operator $L$ 
and the reaction term $R$ admit the form

\begin{equation}
 L \qty(q_c,q_r) = \left(\begin{array}{c} \qty(c+ D_1 \laplacian) q_c \\ D_2 \laplacian q_r\end{array}\right)
\end{equation}
and

\begin{equation}
 R \qty(q_c,q_r) = \left(\begin{array}{c} -a_1 q_c^{\gamma} -a_2
                     q_c^{\beta_c} q_r^{\beta_r} \\ a_1 q_c^{\gamma} +a_2
                     q_c^{\beta_c} q_r^{\beta_r} -d q_r^{\zeta} +
                           B \end{array}\right) .
\end{equation}
For the spatial discretisation a Fourier expansion is applied

\begin{equation}
\label{eq:Pseudospektral-Verfahren-Ansatz}
y \qty(t,x) = \sum_n \varphi \qty(t) \exp(i k_n x) .
\end{equation}
Thus we obtain a system of ordinary differential equations

\begin{equation}
\label{eq:Pseudospektral-Verfahren-ODE}
\dv{\varphi_n}{t} \qty(t) = l_n \varphi_n \qty(t) + R_n \qty(t), \qquad n =
-\frac{N}{2},\dots, \frac{N}{2} ,
\end{equation}
 where $l_n$ is the matrix 

\begin{equation}
l_n = 
\left( \begin{array}{rr}
c- D_1 k_n^2 & 0  \\
0 & -D_2 k_n^2 \\
\end{array} \right) .
\end{equation}
The nonlinear operator $R$ can be expressed with Fourier modes

\begin{equation}
 R \qty(q_c \qty(t),q_r \qty(t)) \approx \sum_{n=-\frac{N}{2}}^{\frac{N}{2}} R_n \qty(t) \exp(i k_n x) .
\end{equation}
Therefore,  the formulation of equation
\eqref{eq:Pseudospektral-Verfahren-ODE} requires a Fourier
transform. In addition after each time step the solution in the
Fourier space, given through \eqref{eq:Pseudospektral-Verfahren-ODE},
can be transformed back and the reaction term can be computed for that
time step. The transformations can be done with a Fast Fourier
Transform.

The system \eqref{eq:Pseudospektral-Verfahren-ODE} can be solved with the 
exponential integrator scheme (\cite{hochbruck_ostermann2010}) of second order 
(ETD2 scheme).
The ETD2 scheme is a two step method. The first step can be computed with 
the exponential integrator scheme of first order 
(ETD1 scheme).

\end{appendix}

\end{document}